\documentclass[11pt,letterpaper,reqno]{article}
\usepackage	[utf8x]{inputenc}
\usepackage [english]{babel}
\usepackage	{amsmath,amssymb,amsthm,mathabx}
\usepackage	{graphicx}
\usepackage	{xcolor}
\usepackage{ulem}

\usepackage[letterpaper]{geometry}
\usepackage{hyperref}
\usepackage{enumitem}
\usepackage{comment}

\DeclareMathAlphabet{\mathpzc}{OT1}{pzc}{m}{it}

\newcommand{\R}{\mathbb{R}}
\newcommand{\T}{\mathbb{T}}

\newcommand{\Z}{\mathbb{Z}}
\newcommand{\N}{\mathbb{N}}

\newcommand{\im}{\mathrm{i}}

\let\sprod\relax\newcommand{\sprod}[2]{\left\langle#1,#2\right\rangle} 

\newtheorem{thm}{Theorem}[section]
\newtheorem{lem}[thm]{Lemma}
\newtheorem{pro}[thm]{Proposition}
\newtheorem{alg}[thm]{Algorithm}

\theoremstyle{definition}
\newtheorem{dfn}[thm]{Definition}
\newtheorem{test}{Test}

\theoremstyle{remark}
\newtheorem{rmk}{Remark}[section]

\title%
{Numerical computation of high-order expansions of invariant manifolds of high-dimensional tori}
\author{Joan Gimeno\footnote{ Departament de Matem\`atiques i Inform\`atica, Universitat de Barcelona,
  Gran Via de les Corts Catalanes, 585, 08007 Barcelona, Spain, \tt joan@maia.ub.es}
\and %
\`Angel Jorba\footnote{ Departament de Matem\`atiques i Inform\`atica, Universitat de Barcelona,
  Gran Via de les Corts Catalanes, 585, 08007 Barcelona, Spain, \tt angel@maia.ub.es}
\and %
Bego\~na Nicol\'as\footnote{ Departament de Matem\`atiques i Inform\`atica, Universitat de Barcelona,
  Gran Via de les Corts Catalanes, 585, 08007 Barcelona, Spain, \tt bego@maia.ub.es (\textit{corresponding author})}
\and %
Estrella Olmedo\footnote{ Barcelona Expert Center (BEC), Institute of Marine Sciences (ICM) and Consejo Superior de Investigaciones Cient\'ificas (CSIC), P. Mar\'itim de la Barceloneta, 37-49, 08003 Barcelona, Spain, \tt olmedo@icm.csic.es}}

\date{\today}

\begin{document}

\maketitle

\begin{abstract}
  In this paper we present a procedure to compute reducible
  invariant tori and their stable and unstable manifolds in stroboscopic Poincar\'e maps.
  The method has two steps. In the first
  step we compute, by means of a quadratically convergent scheme, the
  Fourier series of the torus, its Floquet transformation, and its
  Floquet matrix. If the torus has stable and/or unstable directions,
  in the second step we compute the Taylor-Fourier expansions of the
  corresponding invariant manifolds up to a given order. The paper also discusses
  the case in which the torus is highly unstable so that a multiple shooting
  strategy is needed to compute the torus.
  
  If the order of the Taylor expansion of the manifolds is fixed and $N$ is the
  number of Fourier modes, the whole computational effort (torus and manifolds)
  increases as $\mathcal{O}(N\log N)$ and the memory required behaves as $\mathcal{O}(N)$.
  This makes the algorithm very suitable to compute high-dimensional
  tori for which a huge number of Fourier modes are needed. Besides, the
  algorithm has a very high degree of parallelism. The paper
  includes examples where we compute invariant tori (of dimensions up to 5)
  of quasi-periodically forced ODEs. The computations are run in a
  parallel computer and its efficiency with respect to the number of
  processors is also discussed.
\end{abstract}

\section{Introduction}\label{sec:intro}
The understanding of the dynamics of a given system is usually based
on the knowledge of a subset of solutions that serve as a skeleton for
the global properties. This includes equilibrium points, periodic and
quasi-periodic solutions, and their invariant manifolds.

On the other hand, the accurate description of realistic phenomena
usually requires more sophisticated models. In several cases, some
effects are introduced in the mathematical model as a combination of different periodic
time-dependent perturbations. For example, the description of the
motion of a spacecraft in the neighborhood of the Earth with some level of accuracy requires the
introduction of several effects to the widely used Earth-Moon Circular Restricted Three Body Problem: the effect of the Sun
can be introduced as a periodic time-dependence
(\cite{Huang60,Cronin1964,Andreu98}), and adding some other effects
(like the effects coming from the non-circular motion of the Moon)
results in a quasi-periodic time-dependence. As examples, in
\cite{ESA1, ESA2, ESA3, ESA4}, the authors examine the effect of
perturbing bodies as the sum of periodic perturbations with
incommensurable periods, or in \cite{GomezMasdemontMondelo2002}, where
very detailed models are presented for the motion of a particle in the
Solar system and in the Earth-Moon system, taking into account up to
five natural frequencies of the latter system.
In these cases, we usually have a differential equation,
\begin{equation}
    \dot{x} = f(x,t),
    \label{eq:difeq}
\end{equation}
where $x\in\Omega=\mathring{\Omega}\subset\R ^n$ and $f\colon \Omega \times \R
\rightarrow \R^n$ is a smooth function that we assume to depend
quasi-periodically on $t$. This means that there exists a function
$F\colon \Omega\times \T^{d+1}\rightarrow\R^n$, with $d\geq 1$, such that
\begin{equation*}
    f(x,t) = F(x,\theta_0,\dotsc ,\theta_d),\qquad\theta_i=\omega_i t,
    \quad i=0,1,\ldots,d.
\end{equation*}
Then (\ref{eq:difeq}) can be rewritten as
\begin{equation}
\begin{cases}
\dot{x}=F(x,\theta),\\
\dot{\theta}=\omega,
\end{cases}
\label{eq:DS_d1}
\end{equation}
where $\theta= (\theta_0,\dotsc ,\theta_d)\in\T^{d+1}$ is an angular variable vector and
$\omega=(\omega_0,\ldots,\omega_d)$ is a frequency vector, of
dimension $d+1$, whose components are considered to be linearly
independent over the rationals.  Since $F$ is periodic in each angular variable,
a temporal Poincar\'e section can be defined using one of the angles
of the system, for example $\theta_0=0 \mod{2\pi}$. The associated Poincar\'e
map is the flow of the differential equation from initial time at $0$ to the final time at time
$\delta=\frac{2\pi}{\omega_0}$. Then, redefining the vector $\theta$ as
$(\theta_1,\ldots,\theta_d)$ we can write the Poincar\'e map as
\begin{equation}
\begin{cases}
\bar{x}=P(x,\theta),\\
\bar{\theta}=\theta+\rho,
\end{cases}
\label{eq:DS}
\end{equation}
where $\rho$ is a $d$-dimensional vector with components $\rho_i=\frac{2\pi\omega_i}{\omega_0}$ for 
$i=1,\dotsc ,d$.
The upper bar,~$\;\bar{ }\;$, denotes the image under the Poincar\'e map $P$.

For the discrete system~(\ref{eq:DS}), the simplest possible invariant
sets are invariant tori of dimension $d$ parametrized by the angle
$\theta$, and with frequencies $\rho$. Each of these tori can be
represented by a smooth injective map $\varphi\colon \T^d \rightarrow \R^n$,
which satisfies the following invariance condition
\begin{equation}
P(\varphi(\theta), \theta)=\varphi(\theta+\rho), \qquad \forall\,
\theta \in \T^d.
\label{eq:InvCond}
\end{equation}
Here we are interested in these tori and in an efficient computation of them. Note that, as the dimension of
the torus is higher, the computational effort increases. If the
torus is reducible (see Section~\ref{sec:RedSys}), there are very
efficient numerical methods that allow to compute the torus jointly
with the Floquet change of variables and the reduced Floquet matrix. One of these methods is presented in
\cite{JorbaOlmedo2009}, and it is based on the proofs made in
\cite{JorbaSimo1996}. It is remarkable that the number of operations
and storage requirements are proportional to $N\log N$ and $N$
respectively, where $N$ denotes the number of Fourier modes used to
represent the torus. A side benefit is that the method has a high
degree of parallelism so it can take advantage of modern computers.

In the present work, we extend the methods in \cite{JorbaOlmedo2009} so
that now: {\it i)} the method provides a high-order parametrization of the stable and/or unstable manifolds of the tori if they exist; {\it ii)} we implement a C code that runs in a computer with several processors between which the computations are done concurrently using OpenMP \cite{Dagum1998openmp} instead of the PVM library \cite{GeistBDJMS1994} used in \cite{JorbaOlmedo2009}.

On the one hand, we implement a parametrization
method for computing the Taylor-Fourier expansion of the un/stable invariant manifolds of invariant tori. To implement this
parametrization method on a Poincar\'e map we need to estimate the high-order
derivatives of this Poincar\'e map, which
requires to compute high-order derivatives of the flow of the ODE. This may be very tedious and tricky if we have to compute them by
hand. In particular, since here we focus on the analysis of invariant tori
of maps,  we follow and apply the ideas developed in
\cite{TRAJE2019}, where authors explain how to numerically integrate
high-order derivatives of maps using automatic differentiation with
respect to initial data; a technique called ``jet transport''
\cite{AlessiFJSV08,AlessiFJSV09}. Other authors, \cite{KumarRL21}, also used jet transport to compute Taylor expansions of un/stable manifolds of fixed points of temporal Poincar\'e maps.

On the other hand, the use of parallel computing allows to reduce the
total computational time by using several processors. The parallelization presented in
\cite{JorbaOlmedo2009} was implemented on a cluster of PCs (with a distributed memory), making use
of the PVM library for the communications through an Ethernet network
using a master-slave scheme. Therefore, there was always a penalty time to pay
for those communications. For this reason, the reduction in temporal
costs of the parallelization was relevant for a small number of
processors, but at some point, this reduction stagnates. Nowadays, the
processors involved in our parallelizations belong to the same
computer (hence, with a shared memory), so there is not a communication penalty between the
threads. Both, the computation of the tori and their invariant manifolds are implemented in parallel.

There are many other papers in the literature devoted to the computation of
invariant tori and their manifolds. The computation of a
quasi-periodic solution (with two basic frequencies) of a flow is done, for instance, 
in \cite{DiezJS91} by means of a collocation method.  The
approximation of an invariant tori of a flow is also considered in
\cite{DieciLR91,DieciL95,HuangKM97} using a PDE approach. The computation of 1D
invariant tori of maps is discussed in
\cite{Simo98,CastellaJorba2000}, and an algorithm to obtain their
linear stability is contained in \cite{Jorba2001}. For algorithms
taking advantage of some properties of the linearized normal
behaviour, see \cite{HaroLl05c,HuguetLlS12,CanadellH17,ZhangL18,KumarRL22}.
The main differences between these papers and the present one is that
here we focus on higher dimensional tori of flows, that require a
large amount of computations (and hence, to take advantage of parallel computers) and the use of jet transport to compute
high-order derivatives of Poincar\'e maps.

The structure of this paper is as follows. In
Section~\ref{sec:RedSys}, the reducibility of a quasi-periodic
solution is explained along with a brief summary of the process
developed in \cite{JorbaOlmedo2009} for computing invariant tori for
maps and for obtaining their stability information at the same time. In
Section~\ref{sec:Param}, we introduce the parametrization of hyperbolic
invariant manifolds of those tori, that will be combined with multiple
shooting methods, explained in Section~\ref{sec:MultShoot}, in order
to improve the accuracy of the computations, specially when the
instability of the invariant objects is strong. Details concerning to
the implementation of the computations are included in
Section~\ref{sec:CompImp}.  Section~\ref{sec:Apps} is
devoted to two different applications of the presented technique.
And finally, Section~\ref{sec:conclu} presents some conclusions and future work.

\section{Reducibility of the system}\label{sec:RedSys}
We can analyze the linear behaviour around a torus satisfying \eqref{eq:InvCond}
by taking a small displacement $h\in\R^n$ from a given point $\varphi(\theta)$
of the torus. Applying the map,
\begin{equation*}
P(\varphi(\theta)+h,\theta)=P(\varphi(\theta), \theta)+D_x
P(\varphi(\theta), \theta)h+ \mathcal{O}(\|h\|^2).
\end{equation*}
Renaming this small displacement $h$ as $x$
and defining $A(\theta)=D_x P(\varphi(\theta),\theta)$,
this linear behaviour can be expressed by
the linear skew-product
\begin{equation}
\left. \begin{array}{ll}
\bar{x}=A(\theta)x,\\
\bar{\theta}=\theta+\rho.
\end{array} \right \rbrace
\label{eq:SP}
\end{equation}

\begin{dfn}[Reducible skew product]\label{dfn:reducible}
The system~(\ref{eq:SP}) is said to be reducible if, and only
if, there exists a continuous change of variables  $x=C(\theta)y$ such
that (\ref{eq:SP}) becomes
\begin{equation}
\left. \begin{array}{ll}
\bar{y}=By,\\
\bar{\theta}=\theta+\rho, 
\end{array} \right \rbrace
\label{eq:SP_red}
\end{equation}
where the matrix
\begin{equation}
B=C^{-1}(\theta+\rho)A(\theta)C(\theta),
\label{eq:B}
\end{equation}
does not depend on $\theta$. The matrix $B$ is the {\sl Floquet
matrix} and $x=C(\theta)y$ is the {\sl Floquet transformation}.
\end{dfn}

\begin{dfn}[Discrete reducible torus]
A torus $\varphi$ is reducible when the linear behavior around the torus is reducible.
\end{dfn}

Notice that the Floquet matrix $B$ is a constant matrix that contains
the dynamical information of the system in (\ref{eq:SP}), and then,
its eigenvalues provide the linear stability around the torus. 

\begin{rmk}[Continuous reducible torus]\label{rmk:b_t}
Analogously, $\Psi$ being a quasi-periodic solution, of $d+1$
dimensions, of the system (\ref{eq:difeq}), let us consider the
linearization of the system around the solution as
\begin{equation}
    \dot{x} = \mathpzc{a}(\Psi,t) x
    \label{eq:linear_t}
\end{equation}
where $\mathpzc{a}(\Psi,t)=D_x f (\Psi(t),t)$. We will say that $\Psi$
is a reducible quasi-periodic solution of (\ref{eq:difeq}) if, and only if, there
exists a quasi-periodic change of variables $x=\mathpzc{c}(t)y$ that
transforms (\ref{eq:linear_t}) into
\begin{equation}
    \dot{y} = \mathpzc{b}y,
    \label{eq:linear_t_red}
\end{equation}
where $\mathpzc{b}$ is a constant matrix.

Then, Floquet matrix $B$ for transforming the system (\ref{eq:SP})
into (\ref{eq:SP_red}) corresponds to
\begin{equation*}
    B=\exp \{ \delta \mathpzc{b} \}
    \label{eq:relB_b_t}
\end{equation*}
where $\delta$ denotes the time used to define the temporal Poincar\'e map $P$ and
$\mathpzc{b}$ is the Floquet matrix resulting from transforming
(\ref{eq:linear_t}) into (\ref{eq:linear_t_red}).
\end{rmk}

\begin{rmk}
There is some freedom in choosing the form of matrix $B$. For instance, it can be chosen to be in diagonal (or Jordan) form by composing $C(\theta)$ with a suitable linear (may be complex) transformation, or to be block diagonal to keep $B$ (and $C(\theta)$) real.
\end{rmk}

Although not all the linear skew products (or linear quasi-periodic equations)
are reducible, there are some classical conditions that imply reducibility
\cite{BogoljubovMS76}. It is also remarkable that, under generic conditions,
most of the linear systems that are close to a reducible one are also reducible
\cite{JorbaS92,JorbaRV97}. This applies to the examples considered here.

Assuming the existence of a torus for system~(\ref{eq:DS}) that
satisfies invariance condition (\ref{eq:InvCond}) and that is
reducible, \cite{JorbaOlmedo2009} developed an
iterative method based on Newton iteration (quadratically convergent)
for finding the torus and the Floquet change at the same time. For
this, it is necessary to know suitable seeds for
$\varphi(\theta)$ and $C(\theta)$ (namely $x_0(\theta)$ and $C_0(\theta)$ respectively),
such that $C_0^{-1}(\theta+\rho)D_x P(x_0(\theta),\theta)C_0(\theta)$
is close to a constant matrix $B_0$. Let us suppose by now that these
approximations are chosen. Then, residual magnitudes $y_0(\theta)$
and $Q_0(\theta)$ indicate the error of these
approximations to the real solution, i.e., 
\begin{align}
y_0(\theta) &= x_0(\theta+\rho) - P (x_0(\theta),\theta),
\label{eq:y0} \\
Q_0(\theta) &= C_0^{-1}(\theta+\rho)D_x P(x_0(\theta),\theta)C_0(\theta)-B_0.
\label{eq:Q0}
\end{align}
Then, the norm of these magnitudes is a small quantity of order say 
$\varepsilon$. Let us take, for example, infinite norm,
$\|y_0\|_{\infty}\approx \varepsilon$ and $\|Q_0\|_{\infty}\approx
\varepsilon$. 

The idea is to use the reducibility assumption for finding a better approximation
of the invariant torus, and with it, to improve the Floquet change,
iteratively until the precision of both parametrizations is good
enough.
As all details are carefully described in
\cite{JorbaOlmedo2009}, a very brief description of the procedure is given
here; we address the interested reader to that reference.

Note that in the present work we assume that the torus exists
and that it is real analytic. In \cite{JorbaSimo1996} the proof of its existence is given for
flows and it can be easily translated for maps, as long as some
hypothesis are satisfied. First hypothesis is referred to the
smoothness of the map.  Second hypothesis is that some Diophantine
conditions, involving
the frequency vector $\rho$ and the
eigenvalues $\lambda_1,\dotsc ,\lambda_n$ of the matrix $B_0$
are satisfied. Concretely, it is assumed that there exist real constants $c>0$ and $\gamma >
d-1$ such that
\begin{align}
|\exp( \sprod{\kappa}{\rho} \im) - \lambda_j| &> \frac{c}{|\kappa|^{\gamma}}, \qquad \forall \kappa \in \Z^d \setminus \{0\}, \quad j=1,\dotsc ,n,
\label{eq:Dioph_1} \\
\left|\exp(\langle \kappa, \rho \rangle \im) - \frac{\lambda_j}{\lambda_l}\right| &> \frac{c}{|\kappa|^{\gamma}}, \qquad \forall \, \kappa \in \Z^d \setminus \{0\}, \quad j,l=1,\dotsc ,n,
\label{eq:Dioph_2}
\end{align}
where $\im$ denotes the complex unit, $\sprod{\cdot}{\cdot}$ the
standard scalar product, and $|\kappa| = |\kappa _1| + \dotsb +
|\kappa _d|$. Note that condition (\ref{eq:Dioph_1}) is satisfied if all the
eigenvalues $\lambda_j$ have modulus different from 1, and condition
(\ref{eq:Dioph_2}) is satisfied if all eigenvalues have different
modulus. As at every step of the iterative procedure, the Floquet matrix
changes, a non-degeneracy condition on the eigenvalues is
needed. For this, it is common to make eigenvalues depending on
parameters such that Diophantine condition holds at every
iterative step. In practice, we do not need to verify this condition
but what we have to do is to check if the left-hand side of \eqref{eq:Dioph_1} and
\eqref{eq:Dioph_2} is small. This can be done indirectly by checking the size
of the Fourier modes of the correction given by the Newton method is not too big.

Assuming all these hypothesis to hold and that $x_0(\theta)$ and $C_0(\theta)$
are available (with $\|y_0\|_{\infty}\approx \varepsilon$ and
$\|Q_0\|_{\infty}\approx \varepsilon$), let us summarize the iterative
scheme in \cite{JorbaOlmedo2009} to find good approximations of the torus and of the Floquet change. The iterative scheme is divided in two steps. The first one focuses on computing a better approximation of the torus, and the second one, on
improving the Floquet change and Floquet matrix.

\begin{alg}[Computation of invariant torus, Floquet change, and Floquet matrix] 
\label{alg.torus-floquets} \
\begin{enumerate}
    \item [$\star$] \texttt{Input:} Discrete system as \eqref{eq:DS},
      initial guesses $x _0 (\theta)$, $C _0 (\theta)$, and $B _0$.
    \item [$\star$] \texttt{Output:} Torus $\varphi(\theta)$, Floquet
      change $C(\theta)$, and Floquet matrix $B$.
\end{enumerate}
First step:
\begin{enumerate}
\renewcommand{\theenumi}{1.\roman{enumi}}
    \item Compute the error $y _0(\theta) = x _0(\theta + \rho) - P(x
      _0(\theta), \theta)$.
    \item Compute the function $g(\theta)=-C_0^{-1}(\theta+\rho)
      y_0(\theta)$.
    \item Find $u$ that verifies $u(\theta+\rho)=B_0
      u(\theta) + g(\theta)$. For this, we expand functions $g$
      and $u$ in real Fourier series (the expansion can be done
      in complex Fourier series, but we work with real expansions in
      the computer programs):
    \begin{equation*}
        \begin{aligned}
        g(\theta) &=\frac{g^{(0)}}{2} + \sum\limits_{\kappa\neq 0} g_{\kappa}^{(c)} \cos \langle \kappa,\theta \rangle + g_{\kappa}^{(s)} \sin \langle \kappa,\theta \rangle, \\
        u(\theta) &=\frac{u^{(0)}}{2} + \sum_{\kappa\neq 0} u_{\kappa}^{(c)} \cos \langle \kappa,\theta \rangle + u_{\kappa}^{(s)} \sin \langle \kappa,\theta \rangle,
        \end{aligned}
        \label{eq:g_varphi}
    \end{equation*}
    where $\kappa \in \N ^d$ and $\langle \kappa,\theta \rangle=\kappa_1 \theta_1+ \dotsb
    +\kappa_d \theta_d $, and solve the following system to find
    Fourier coefficients $u^{(0)}$, $u_{\kappa}^{(c)}$,
    and $u_{\kappa}^{(s)}$
    \begin{equation}
        \begin{aligned}
        (Id-B_0)\frac{u^{(0)}}{2} &= \frac{g^{(0)}}{2},\\
        (B_0^2 -2 \cos \langle \kappa,\rho \rangle B_0 + Id)u_{\kappa}^{(c)} &= (\cos \langle \kappa,\rho \rangle Id + B_0)g_{\kappa}^{(c)} - \sin \langle \kappa,\rho \rangle g_{\kappa}^{(s)},\\
        (B_0^2 -2 \cos \langle \kappa,\rho \rangle B_0 + Id)u_{\kappa}^{(s)} &= ( \cos \langle \kappa,\rho \rangle Id + B_0)g_{\kappa}^{(s)} + \sin \langle \kappa,\rho \rangle g_{\kappa}^{(c)},
        \end{aligned}
        \label{eq:coho_torus}
    \end{equation}
    where $Id$ denotes the identity matrix. 
    \item Compute $h(\theta)=C_0(\theta) u(\theta)$.
    \item Compute $x_1(\theta)=x_0(\theta)+h(\theta)$, that is the new
      approximation of the torus, such that $\|y_1\|_{\infty}\approx
      \varepsilon^2$, with $y_1$ defined like in (\ref{eq:y0}).
\end{enumerate}

\noindent Second step:
\begin{enumerate}
\renewcommand{\theenumi}{2.\roman{enumi}}
    \item Compute the matrices $R(\theta)=C_0^{-1}(\theta +\rho) D_x
      f(x_1(\theta), \theta)C_0(\theta) - B_0$,
      $\tilde{R}(\theta)=R(\theta) - \textrm{Avg}(R)$, and
      $B_1=B_0+\textrm{Avg}(R)$, where $\textrm{Avg}(R)$ is the
      average of the map $\theta \mapsto R(\theta)$, that is,
    \begin{equation*}
        \textrm{Avg}(R)=\frac{1}{(2\pi)^d} \int_{\T^d} R(\theta) \, d\theta
    \end{equation*}
    and $B_1$ is the new approximation to the Floquet matrix.
    \item Find the matrix valued function $H$ that verifies
      $H(\theta+\rho)B_1 -B_1 H(\theta)=\tilde{R}(\theta)$. For this,
      we expand $R$ and $H$ in real Fourier series:
    \begin{equation*}
        \begin{aligned}
        H(\theta) &= \sum\limits_{\kappa \neq 0} H_{\kappa}^{(c)} \cos \langle \kappa,\theta \rangle + H_{\kappa}^{(s)} \sin \langle \kappa,\theta \rangle, \\
        \tilde R(\theta) &= \sum_{\kappa\neq 0} R_{\kappa}^{(c)} \cos \langle \kappa,\theta \rangle + R_{\kappa}^{(s)} \sin \langle \kappa,\theta \rangle,
        \end{aligned}
        \label{eq:H_varphi}
    \end{equation*}
    and we solve the following system to find Fourier coefficients
    $H_{\kappa}^{(c)}$ and $H_{\kappa}^{(s)}$ (note that $H^{(0)}=0$)
    \begin{equation}
        \begin{aligned}
        (H_{\kappa}^{(c)} \cos \langle \kappa,\rho \rangle + H_{\kappa}^{(s)} \sin \langle \kappa,\rho \rangle) B_1 -B_1 H_{\kappa}^{(c)}&= R_{\kappa}^{(c)},\\
        (H_{\kappa}^{(s)} \cos \langle \kappa,\rho \rangle - H_{\kappa}^{(c)} \sin \langle \kappa,\rho \rangle) B_1 -B_1 H_{\kappa}^{(s)}&= R_{\kappa}^{(s)}.
        \end{aligned}
        \label{eq:coho_floquet}
    \end{equation}
    \item Compute $C_1(\theta)=C_0(\theta) (Id + H(\theta))$, that is
      the new approximation of the Floquet transformation, such that
      $\|Q_1\|_{\infty}\approx \varepsilon^2$, with $Q_1$ defined like
      in (\ref{eq:Q0}).
\end{enumerate}
Once we have $x_1(\theta)$, $B_1$, and $C_1(\theta)$, we keep on iterating
until either the norms of $y$ and $Q$ are small enough or the differences between one step and the previous one are small enough.
\end{alg}

Note that the computation of each pair of coefficients
($u_{\kappa}^{(c)}$, $u_{\kappa}^{(s)}$) 
in
(\ref{eq:coho_torus}) is independent for each $\kappa$. The same
happens for each pair ($H_{\kappa}^{(c)}$, $H_{\kappa}^{(s)}$) in
(\ref{eq:coho_floquet}). That makes these linear systems very suitable
for their computational resolution in parallel. With this, the
dimension of each linear system depends on the dimension of the phase
space and the number of linear systems to be solved on the number of
Fourier modes used. On the other hand, the evaluation of the map $P$ and $D _x P$ can be performed independently of each $\theta$, which leads to a straightforward parallelization.

\section{Parametrization of invariant manifolds} \label{sec:Param}
The quasi-periodic solutions of interest in this work are those that
have a saddle part, i.e., those that have un/stable invariant
manifolds associated. In that case, the eigenfunctions and eigenvalues associated with their hyperbolic directions are typically used to compute linear approximations to the invariant manifolds. Here, we present a
procedure to compute a high-order approximation of the stable or unstable invariant manifolds for a torus in the map \eqref{eq:DS}. 
In fact, the procedure presented here can be seen as an extension of the one in \cite{JorbaNicolas2021} for the computation of the un/stable invariant manifolds of invariant curves.

Let us write the invariant manifold of an invariant torus as a formal
Taylor-Fourier expansion in terms of two parameters, a parameter $\sigma\in \R$ and the
angle vector $\theta \in \T^d$:
\begin{equation}
W(\theta, \sigma) = a_0(\theta) + a_1(\theta) \sigma + \sum_{k\geq 2} a_k(\theta) \sigma^k,
\label{eq:param}
\end{equation}
where $a_k$ is a function from $\T^d$ to $\R^n$.

Here we consider the case in which $\sigma$ is one-dimensional. 
The case of invariant
manifolds with several hyperbolic directions
requires a parameter $\sigma _j$ for each direction. In such a case, $k$ in \eqref{eq:param} becomes a multi-index and the Taylor-Fourier polynomial manipulations have to be modified accordingly. The case of invariant
manifolds with several hyperbolic directions associated with fixed
points of Poincar\'e maps is detailed in \cite{TRAJE2019}.

In the particular case that the system is Hamiltonian, the stable eigenvalue
corresponds to the inverse of the unstable one. And, since many of the
classical mechanical systems show a symmetry when inverting the
time, sometimes it is possible to have both parametrizations,
the stable and the unstable, by just computing one of them and applying the
corresponding symmetry. Here, we are going to explain the general case, valid for Hamiltonian and non-Hamiltonian systems and without considering, in advance, any symmetry.

\subsection{Computation of the invariant manifold}

The invariance condition for an invariant manifold $W$ of a torus of a Poincar\'e map $P$ like \eqref{eq:DS} can be written as
\begin{equation}
P(W(\theta,\sigma), \theta)=W(\theta+\rho, \lambda \sigma),
\label{eq:InvEq_manitorus}
\end{equation}
where $\rho$ is the frequency vector and $\lambda$ is a real hyperbolic
eigenvalue. The right hand side of \eqref{eq:InvEq_manitorus} can be written as
\[
 W(\theta+\rho, \lambda \sigma)= a_0(\theta + \rho) + a_1(\theta +
 \rho) \lambda \sigma + \sum_{k\geq 2}a_k(\theta+\rho) (\lambda
 \sigma)^k.
\]
Each function $a_k$ will be obtained by solving the invariance
condition~(\ref{eq:InvEq_manitorus}) order by order. Let us denote by
$W_m(\theta,\sigma)$ the truncated Taylor-Fourier series up to order
$m$ in a neighbourhood of $\sigma=0$,
\begin{equation}
    W_m(\theta,\sigma)=a_0(\theta)+a_1(\theta)\sigma+ \sum_{k=2}^m
    a_k(\theta) \sigma^k.
\label{eq:param_m}
\end{equation}
 
From equations (\ref{eq:param}) and
(\ref{eq:InvEq_manitorus}), we see that $a_0$ is the
parametrization of the invariant torus and $a_1$ is the hyperbolic
eigenfunction, that are obtained as truncated real Fourier series
following the iterative scheme in the previous section, see also Algorithm~\ref{alg.mani-unstab}. More precisely, we have obtained the eigenvectors of the Floquet matrix $B$. Then, if $v$  is the eigenvector of $B$ associated with $\lambda$, the eigenfunction $a_1$ of the torus is $C(\theta) v$.
The torus and its eigendirection 
give the linear approximation to the invariant manifolds, which is the expansion in (\ref{eq:param}) up to
order one.

Assume that we know the parametrization $W_{m-1}$ up to order $m-1$,
i.e. the functions $a_0,\dotsc ,a_{m-1}$ are known. We want to
find the function $a_m$ involved in $W_m$. For this, we apply the
Poincar\'e map to $W_m (\theta, \sigma) = W _{m-1}(\theta , \sigma) + a _m(\theta) \sigma ^m$,
\begin{equation}
\begin{split}
P(W_m(\theta,\sigma), \theta) &= P(W_{m-1}(\theta,\sigma), \theta) + D_x P(W_{m-1} (\theta, \sigma), \theta) a_m(\theta) \sigma^m + \mathcal{O} (\sigma^{m+1}) \\
 &= P(W_{m-1}(\theta,\sigma), \theta) + D_x P(a_0(\theta), \theta) a_m(\theta) \sigma^m + \mathcal{O} (\sigma^{m+1}).
 \label{eq:fWm}
\end{split}
\end{equation}
The Taylor-Fourier expansion of the invariant manifold $W_{m-1}$ under
the Poincar\'e map is
\[
 P(W_{m-1}(\theta,\sigma), \theta) = W_{m-1}(\theta + \rho, \lambda
 \sigma) + b_m (\theta) \sigma^m + \mathcal{O} (\sigma^{m+1}),
\]
where $b_m$ is a $\theta$-dependent function at order $m$ that comes
from the evaluation of the invariant manifold at order $m-1$. For the
computation of this function we make use of the jet transport
technique, detailed in Section~\ref{subsec:JetTransport}.
Note that we only need to compute the expansion of $P(W_{m-1}(\theta,\sigma), \theta)$ w.r.t. $\sigma$ up to order $m$.

Now we insert the last expression in (\ref{eq:fWm}) and impose the invariance
condition (\ref{eq:InvEq_manitorus}) up to order $m$:
\begin{equation*}
    b_m(\theta) \sigma^m + D_x P (a_0(\theta), \theta) a_m (\theta)
    \sigma^m = a_m(\theta+\rho)\lambda^m \sigma^m.
    \label{eq:am}
\end{equation*}
If the number of harmonics used is low, functions $a_k$ can be computed
by solving the above linear system in which the matrix $D_x
P(a_0(\theta),\theta)$ is involved.
If the dimension of the torus is high, so is the number of harmonics and then solving this system directly is not feasible.
However, when the system is reducible we can express the previous equation in terms
of the Floquet matrix $B$, such that the final system to solve offers
again a high degree of parallelism.

Therefore, we introduce the Floquet change 
$a_m(\theta) = C(\theta) u_m(\theta)$,
\[
 b_m(\theta) + D_x P (a_0(\theta), \theta) C(\theta) u_m (\theta) =
 C(\theta +\rho) u_m(\theta + \rho) \lambda^m,
\]
and multiplying by $C^{-1}(\theta+\rho)$ on both sides, it leads to
\begin{equation}
    C^{-1}(\theta+\rho) b_m(\theta) + B u_m (\theta) = u_m(\theta +
    \rho) \lambda^m.
\label{eq:um}
\end{equation}
Under the generic condition of non-resonance, detailed in the Lemma~\ref{lem:ex_v}, \eqref{eq:um} determines uniquely the function $u_m$, that gives $a_m$ through the Floquet transformation.

So, let us assume that $g _m(\theta) = C ^{-1}(\theta + \rho) b _m(\theta)$ admits a (real) Fourier series expansion, that is,
\begin{equation*}
     g_m(\theta) =\frac{g^{(0)}}{2} + \sum\limits_{\kappa \neq 0} g_{\kappa}^{(c)} \cos \langle \kappa,\theta \rangle + g_{\kappa}^{(s)} \sin \langle \kappa,\theta \rangle.
\end{equation*}
Then, we have to find the coefficients of another Fourier expansion
 \begin{equation*}
     u_m(\theta) =\frac{u^{(0)}}{2} + \sum_{\kappa\neq 0} u_{\kappa}^{(c)} \cos \langle \kappa,\theta \rangle + u_{\kappa}^{(s)} \sin \langle \kappa,\theta \rangle ,
 \end{equation*}
such that \eqref{eq:um} is satisfied. Imposing the equation \eqref{eq:um} on the Fourier coefficients and using a Cramer-block method, we end up with the following linear system of cohomological equations depending on $\kappa$,
\begin{equation}
\begin{aligned}
(\lambda^m Id-B)\frac{u^{(0)}}{2} &= \frac{g^{(0)}}{2},\\
(B^2 -2\lambda^m \cos \langle \kappa,\rho \rangle B + \lambda^{2m}Id)u_{\kappa}^{(c)} &= (\lambda^m \cos \langle \kappa,\rho \rangle Id + B)g_{\kappa}^{(c)} - \lambda^m \sin \langle \kappa,\rho \rangle g_{\kappa}^{(s)},\\
(B^2 -2\lambda^m \cos \langle \kappa,\rho \rangle B + \lambda^{2m}Id)u_{\kappa}^{(s)} &= (\lambda^m \cos \langle \kappa,\rho \rangle Id + B)g_{\kappa}^{(s)} + \lambda^m \sin \langle \kappa,\rho \rangle g_{\kappa}^{(c)}.
\end{aligned}
\label{eq:coho_eq}
\end{equation}

The linear systems in \eqref{eq:coho_eq} are solvable as long as $B^{(\kappa)} = B ^2 - 2
\lambda ^m\cos \langle \kappa, \rho \rangle B + \lambda^{2m} Id$ is invertible for all $\kappa$.
If $\mu$ is an eigenvalue of $B$, then $B ^{(\kappa)}$ has eigenvalues of the form
\begin{equation*}
 \mu ^2 - 2 \lambda ^m \cos \langle \kappa, \rho \rangle \mu + \lambda ^{2m}
\end{equation*}
which makes $B^{(\kappa)}$ invertible whenever $\mu$ is different to
$\lambda ^m \exp(\pm \langle \kappa, \rho \rangle \im)$.

We have then proved the following lemma:

\begin{lem}
 Let $B$ be a Floquet matrix associated with the frequency vector $\rho$ on $\T ^d$ and let $|\lambda| \ne 1$ be a real number satisfying that for each eigenvalue $\mu$ of $B$ and a fixed $m \in \N$, $m \geq 2$,
 \begin{equation}
     \mu \ne \lambda ^m\exp(\pm \langle \kappa, \rho \rangle \im), \qquad \forall\, \kappa \in \N ^d.
\label{obs:mus_B}
 \end{equation}
 Then for all smooth function $g _m$ on $\T ^d $,
 there exists a unique smooth function $u _m$ such that 
 \begin{equation}
  \lambda^m u_m(\theta+\rho)=B u_m(\theta) + g_m(\theta).
  \label{eq:tosolve}
 \end{equation}
\label{lem:ex_v}
\end{lem}

\begin{rmk}
 Note that \eqref{obs:mus_B} is always satisfied when $\lambda$ is the dominant eigenvalue of the Floquet matrix.
\end{rmk}

Expressions in (\ref{eq:coho_eq}) recall those in
(\ref{eq:coho_torus}) and (\ref{eq:coho_floquet}). Therefore, it is clear
that the computation of each pair of coefficients 
$(u^{(c)}_{\kappa},u^{(s)}_{\kappa})$ of $u_m$
is
independent to each other. That makes the proposed invariant manifold computation highly parallelizable as it was the Algorithm~\ref{alg.torus-floquets}. Hence, to find the unknowns 
$(u^{(c)}_{\kappa},u^{(s)}_{\kappa})$, 
we solve a large
number of small dimensional linear systems at the same time; the
dimension of each linear system depends on the dimension of the phase
space and the number of systems only depends on the number of Fourier
modes used for the approximation of 
$g _m$ and $u _m$.

Notice that, when we start the computation of the functions $a_k$ for $k
\geq 2$, we use the same number of Fourier modes $N _i$ for each of the angular dimensions $\theta_i$ with $i = 1, \dotsc, d$, as for the torus and the Floquet change Fourier
series. However, it may happen that those numbers of modes, that were enough for discretizing accurately the torus and the Floquet change, may not be enough for discretizing some of the parametrization functions $a_k$ for $k \geq 2$. If this happens for a given order, it is necessary to increase the number of Fourier modes from this order on. This has not been the case in the examples of Section~\ref{sec:Apps}, where we have checked the size of the Fourier modes after the
computation. Varying the number of Fourier modes during the computation has an
extra penalty that depends on the number of Fourier modes added.
An example of a varying number of Fourier modes during the computation of the manifold
can be found in \cite{JorbaNicolas2021}.

Algorithm~\ref{alg.mani-unstab} summarizes the process explained
above.

\begin{alg}[Invariant manifold of a torus through its Floquet transformation] \
\label{alg.mani-unstab}
\begin{enumerate}
\renewcommand{\theenumi}{\arabic{enumi}}
    \item [$\star$] \texttt{Input:} Discrete system as in
      \eqref{eq:DS}, torus $\varphi(\theta)$, Floquet change
      $C(\theta)$, Floquet matrix $B$, real eigenvalue $|\lambda | \ne
      1$ of $B$, and its eigenvector $v$.
    \item [$\star$] \texttt{Output:} Coefficients $a_k(\theta)$ for $k\geq 2$ verifying
      \eqref{eq:InvEq_manitorus}.
    \item $a _0(\theta) \gets \varphi(\theta)$.
    \item $a _1(\theta) \gets C(\theta) v$.
    \item For $k = 2, 3, \dotsc$
    \begin{enumerate}
\renewcommand{\theenumii}{\alph{enumii}}
\renewcommand{\labelenumii}{\theenumii)}
        \item \label{alg.mani-unstab-a} $b _0(\theta) + \dotsb + b _{k}(\theta) \sigma^{k} \gets
          P(a _0(\theta) + \dotsb + a _{k-1}(\theta) \sigma^{k-1} + \mathbf{0}\sigma ^k,
          \theta)$ using jet transport.
        \item $g _k(\theta) \gets C ^{-1}(\theta + \rho) b _k(\theta)$.
        \item Find $u _k(\theta)$ such that $\lambda ^k u _k(\theta +
          \rho) = B u _k(\theta) + g _k(\theta)$ using
          \eqref{eq:coho_eq}.
        \item $a _k(\theta) \gets C(\theta) u _k(\theta)$.
    \end{enumerate}
\end{enumerate}
\end{alg}

\subsection{Stable invariant manifold}
\label{sec:param_stable}
The procedure introduced above is valid for both, the stable and the unstable invariant manifold computation. However, when the Poincar\'e map is
 applied forward in time to compute the stable manifold, it approaches the
 torus and also the unstable invariant manifold. This computation affects the numerical accuracy of the stable manifold by increasing the numerical errors due to the effect of the unstable direction. This effect is more relevant when the unstable direction is strong. Because of that, it is more accurate to obtain the parametrization of the
 stable invariant manifold using the inverse of the Poincar\'e map in \eqref{eq:DS},
\begin{equation}
\begin{cases}
x=P^{-1}(\bar{x},\bar{\theta}),\\
\theta=\bar{\theta}-\rho. 
\end{cases}
\label{eq:InvDS}
\end{equation}

In the case of the stable invariant manifold, i.e. the reals $|\lambda| < 1$ eigenvalues of $B$ in \eqref{eq:SP_red}, the invariance condition
(\ref{eq:InvEq_manitorus}) is written for \eqref{eq:InvDS} as,
\begin{equation}
P^{-1}(W(\theta+\rho,\sigma), \theta + \rho)=W \bigl(\theta,  \frac{\sigma}{\lambda} \bigr).
\label{eq:InvEq_manitorus_inverse}
\end{equation}

We proceed as before, we consider a formal power expansion of $W$ in \eqref{eq:param} and solve (\ref{eq:InvEq_manitorus_inverse}) order by order to obtain the functions $a_k$ that parametrize the stable invariant manifold. 

Notice that in that case, we must introduce a Floquet change that removes the angle dependence when the dynamics is moving backward in time, i.e. when we apply $P^{-1}$. The Floquet transformation for the torus, and the torus itself, in the map $P$ and in the map $P^{-1}$ are related through a phase equal to the vector $\rho$.

Following the Algorithm~\ref{alg.torus-floquets}, the invariant torus and its eigenfunction (order zero and one of the parametrization) are obtained by application of $P$. If we want to use them for the parametrization of the stable manifold, where we apply $P^{-1}$, we have to re-parametrize them as $a_0(\theta+\rho)$ and
$a_1(\theta+\rho)$. Therefore, in this case we look for the functions $a_k$ with $k\geq 2$ shifted a quantity $\rho$.

Then, assuming that we know the parametrization up to order $m-1$, in order to find the function $a_m(\theta) = C(\theta) u_m(\theta)$, we find $u_m(\theta)$ satisfying
\begin{equation}
    C^{-1}(\theta)b^{-}_m(\theta) + B^{-1}u_m (\theta+\rho) = \lambda^{-m}u_m (\theta) ,
\end{equation}
where $b^{-}_m(\theta)$ denotes the term of order $m$
resulting from the evaluation of the invariant manifold up to order
$m-1$ by the inverse Poincar\'e map.
Now, multiplying by $B$ and by $\lambda^{m}$ the last expression, we have
\begin{equation}
 \lambda^m u_m (\theta+\rho) = B u_m (\theta) + g_m(\theta), 
 \label{eq:tosolve_s}
\end{equation}
that has the same form as \eqref{eq:tosolve} with $g_m (\theta) =- \lambda^m B C^{-1}(\theta)b^{-}_m(\theta)$. Therefore, relaying on Lemma~\ref{lem:ex_v} and the solution of the linear systems in \eqref{eq:coho_eq}, there exists the
function $u_m\colon \T^d \mapsto \R^n$, evaluated in $u_m(\theta+\rho)$,
that satisfies (\ref{eq:tosolve_s}). Note that, the condition in \eqref{obs:mus_B} remains the same.

\begin{alg}[Stable invariant manifold of a torus through its Floquet transformation] \
\label{alg.mani-stab}
\begin{enumerate}
\renewcommand{\theenumi}{\arabic{enumi}}
    \item [$\star$] \texttt{Input:} Discrete system as in
      \eqref{eq:DS}, torus $\varphi(\theta)$, Floquet change
      $C(\theta)$, Floquet matrix $B$, real eigenvalue $|\lambda| <
      1$ of $B$, and its eigenvector $v$.
    \item [$\star$] \texttt{Output:} Coefficients $a _k(\theta + \rho)$ for $k \geq 2$ 
      verifying
      \eqref{eq:InvEq_manitorus_inverse}.
    \item $a _0(\theta + \rho) \gets \varphi(\theta + \rho)$.
    \item $a _1(\theta + \rho) \gets C(\theta + \rho) v$.
    \item For $k = 2, 3, \dotsc$
    \begin{enumerate}
\renewcommand{\theenumii}{\alph{enumii}}
\renewcommand{\labelenumii}{\theenumii)}
        \item \label{alg.mani-stab-a} $b _0(\theta) + \dotsb + b _{k}(\theta) \sigma^{k} \gets
          P^{-1}(a _0(\theta + \rho) + \dotsb + a _{k-1}(\theta + \rho) \sigma^{k-1} + \mathbf{0}\sigma ^k,
          \theta +\rho)$ using jet transport.
        \item $g _k(\theta) \gets -\lambda ^k B C ^{-1}(\theta ) b _k(\theta)$.
        \item Find $u _k(\theta)$ such that $\lambda ^k u
          _k(\theta + \rho) = B u _k(\theta) + g _k(\theta)$ using
          \eqref{eq:coho_eq}.
        \item $u _k(\theta + \rho) \gets u _k(\theta )$.  
        \item $a _k(\theta + \rho) \gets C(\theta + \rho) u _k(\theta + \rho)$.
    \end{enumerate}
\end{enumerate}
\end{alg}

\subsection{Scaling factor}
The parametrization of the invariant manifold is not unique. In particular, if $W(\theta, \sigma)$ is a parametrization, then $W(\theta, c \sigma')$ for $c \ne 0$ is also a parametrization
w.r.t. a new parameter $\sigma'$. The role of $c$ is to rescale the coefficients $a _k(\theta)$ in \eqref{eq:param} as $c ^k a _k(\theta)$. To minimize the error propagation we want to avoid the norms of the coefficients to grow or decrease too fast with $k$ \cite{Richardson80,FalcoliniL92}. 

If the radius of convergence w.r.t. $\sigma$ is known, $\varrho$, we can use the scaling $c = \varrho $. If $\varrho$ is not known, we can run the algorithm for the manifold twice; one to estimate it and a second one to rescale. Note that if after the first run the estimated radius of convergence is not far from $1$, we accept the computed expansion without the second run. 
On the other hand, if the radius is far from $1$, we can decide to recompute the manifold with a proper scaling. 

Finally, a simple way to rescale the expansion is to scale the eigenvector of the Floquet matrix to have norm $c$ and run the algorithm.

\section{Multiple shooting}\label{sec:MultShoot}

There are invariant objects so unstable that it is impossible to
integrate accurately the flow around them during the
time involved in the Poincar\'e map. In these cases, it is convenient to
split the time into a certain number of temporal sections, such that
the integration time between each two consecutive sections is considerably
reduced and so the propagation of the numerical error. These methods are commonly known as multiple shooting or parallel shooting methods and they have been widely used to compute highly unstable  periodic orbits \cite{StoerB02}. A version of this methodology
applied to the computation of invariant tori and their Floquet changes
can be found in \cite{OlmedoPhD}.

The section $\Sigma = \{(x,\theta) \in \R ^{n} \times \T ^d \colon \theta _0 \equiv 0 \}$ is the one considered in the map $P$ of \eqref{eq:DS} and in the case that the underlying problem involves an ODE, the map $P \colon \Sigma \to \Sigma$ can be given in terms of a flow map of period $\delta = 2\pi / \omega _0$. 
\newline \newline
In a multiple shooting method, instead of considering only one section $\Sigma$, a finite number of sections are considered such that the last section coincides with the first one. More precisely, let us consider the $r$ sections \[\Sigma _j = \{(x, \theta) \in \R ^{n} \times \T ^d \colon \theta _0 \equiv (j-1) \delta / r \}, \qquad j = 1, \dotsc, r\]
and $\Sigma _{r+1}$ is exactly $\Sigma _1$. For each of two consecutive of these sections, we consider maps $P _j \colon \Sigma _j \to \Sigma _{j+1}$ that, for instance, in an ODE case, they can be defined by the same flow considered in $P$ but stopping at the different values of $\theta _0$ in $\Sigma _j$.
An invariant torus can be represented by a set of smooth
maps $\varphi_j:\T^d\rightarrow\Sigma_j$ which are the intersections of
a $(d+1)$-dimensional torus of the flow with each of the sections.
Now the invariance equation \eqref{eq:InvCond} for the torus splits in conditions for the new (unknown) functions $\varphi _1, \dotsc, \varphi _r$
imposing that $\varphi_j$ is mapped on $\varphi_{j+1}$ by $P_j$
($j=1,\ldots,r-1$) and that $\varphi_r$ is mapped by $P_r$
on $\varphi_1$, but with a shift of $\rho$ in $\theta$.

Our goal is to choose an expression that only changes the evaluation of the
$P$ and $D _x P$ and the other steps in Algorithms
\ref{alg.torus-floquets}, \ref{alg.mani-unstab}, and \ref{alg.mani-stab} remain the same. To keep the same shifting over the
$r$ new sections, we apply shifts to the equations 
\begin{equation*}
    \begin{split}
        P _j (\varphi _j(\theta + (j-1)\rho/r), \theta + (j-1)\rho/r) &= \varphi _{j+1} (\theta + j\rho /r), \qquad j = 1, \dotsc, r-1, \\
        P _r (\varphi _r(\theta+(r-1)\rho/r), \theta+(r-1)\rho/r ) &= \varphi _{1} (\theta + \rho).
    \end{split}
\end{equation*} 
where the interpretation of these $P _j$ in an ODE scenario is a flow that advances certain time in $\theta _0$ to reach the section $\Sigma _{j+1}$. Thus, the invariance equation after the shifts is written as
\begin{equation}
\label{eq.inv-eq-ps}
    \begin{split}
        P _j (\varphi _j(\theta), \theta) &= \varphi _{j+1} (\theta + \rho /r), \qquad j = 1, \dotsc, r-1, \\
        P _r (\varphi _r(\theta), \theta) &= \varphi _{1} (\theta + \rho/r).
    \end{split}
\end{equation}
Note that one can see a multiple shooting approach as a single shooting but with a larger phase space.
The equations \eqref{eq.inv-eq-ps} do not recover $P$ in
\eqref{eq:InvCond} by a direct composition since they must be
alternated with rotation operators. Similarly, the $\varphi$ in
\eqref{eq:InvCond} can be obtained from the $\varphi _1, \dotsc,
\varphi _r$ in \eqref{eq.inv-eq-ps} by undoing rotations. The
Lemma~\ref{lem.PPj} makes explicit all these rotations and, in particular, says that $D _x P$ is a product of differentials of $P _j$ with interlaced rotations.

\begin{lem}
\label{lem.PPj}
 Let $\alpha$ be an angle and let $T _\alpha$ be the operator defined
 as $T _\alpha x(\theta) = x (\theta +\alpha)$. Then the map $P$ at the torus $\varphi$ in
 \eqref{eq:InvCond} and $P _j$ in \eqref{eq.inv-eq-ps} at the torus $\varphi _j$ are related by
\begin{equation*}
 P = T _{\rho - \rho /r} \circ P _r \circ T _{-\rho /r} \circ P _{r-1} \circ \dotsb \circ  T _{-\rho /r} \circ P _{1}.
\end{equation*}
\end{lem}

Some authors like \cite{DuartePhD,RosalesJJC21a} have computed invariant tori for flows that depend on time in a periodic way. In this case, the maps $P_j$ are autonomous (they do not depend on a perturbing angle $\theta$) and then
the parametrizations $\varphi_j$ are not unique. This allows for different
expressions in (\ref{eq.inv-eq-ps}).

\subsection{Reducibility and multiple shooting}
The linear skew-product associated to the linearization around a torus found
with this multiple shooting has a uniform rotation over the sections. Thus,
\eqref{eq.inv-eq-ps} has a linear behaviour expressed by
\begin{equation}
    \label{eq:SPj}
\begin{split}
 \bar{X}_{j+1} &= A _j (\theta) X_{j},   \qquad j=1,\dotsc ,r-1, \\
 \bar{X}_1 &= A _r(\theta) X_{r},  \\
 \bar{\theta} &= \theta + \rho/r,
\end{split}
\end{equation}
with $A _j(\theta) = D_x P_{j} (\varphi_{j} (\theta), \theta)$ for $j = 1, \dotsc, r$.  

Following the Definition~\ref{dfn:reducible}, \eqref{eq:SPj} is
reducible if, and only if, there exist a change of variables of the form
$ X _j = C _j(\theta) Y _j$ for $j=1,\dotsc ,r$ such
that~(\ref{eq:SPj}) becomes
\begin{equation}
    \label{eq:SPj_red}
\begin{split}
  \bar{Y}_{j+1} &= B_j Y_j,   \qquad j=1,\dotsc ,r-1  \\
  \bar{Y}_{1} &= B_r Y_r,  \\
  \bar{\theta} &= \theta + \rho/r,
\end{split}
\end{equation} 
where the matrices $B _1, \dotsc, B _r \in \R ^{n\times n}$ are defined by 
\begin{equation} \label{eq:SPj_B}
\begin{split}
  B _{j} &= C _{j+1}(\theta + \rho /r) ^{-1} A _j(\theta) C _j(\theta),   \qquad j=1,\dotsc ,r-1  \\
  B _{r} &= C _1(\theta + \rho /r) ^{-1} A _r(\theta) C _r(\theta),
\end{split}
\end{equation}
and they do not depend on $\theta$.

\begin{rmk}[Matrix-form] \label{rmk.mat-form-ps}
 The linearization around a torus can also be formulated in a matrix-block form of a higher dimensional problem, that is, $n \cdot r$ dimension. Thus in practice, Algorithm~\ref{alg.torus-floquets} can be used with
 multiple shooting just considering $\rho/r$ instead of $\rho$ and
 evaluating $P _j$ and $D _x P _j$ following
 \eqref{eq.inv-eq-ps}. Indeed, if we consider the block matrices 
 \begin{align}
  \label{eq.mat-form-ps}
  \tilde A(\theta) &= 
  \begin{pmatrix}
   & & & & A _{r} \\
   A _1 \\
   &   \ddots \\
   & & A _{r-2} \\ 
   & & & A _{r-1}
  \end{pmatrix}(\theta), &
  \tilde B &= 
  \begin{pmatrix}
   & & & & B _{r-1} \\
   B _r \\
   & B _{1} \\
   & &  \ddots \\
   & & & B _{r-2} 
  \end{pmatrix},\\ \nonumber
  \tilde C(\theta) &= 
  \begin{pmatrix}
   & C _1 \\
   & &  \ddots \\
   & & & C _{r-2} \\
   & & & & C _{r-1} \\
   C _{r}
  \end{pmatrix}(\theta), & \tilde C ^{-1}(\theta) &= 
  \begin{pmatrix}
   & & & & C _{r}^{-1} \\
   C _1^{-1} \\
   & \ddots \\
   & & C _{r-2}^{-1} \\
   & & & C _{r-1}^{-1}
  \end{pmatrix}(\theta),
 \end{align}
 then using a little bit more memory to keep the zeros for each
 $\theta$ we can directly use Algorithm~\ref{alg.torus-floquets}.
\end{rmk}

As a consequence of the Remark~\ref{rmk.mat-form-ps}, we have the
following straightforward lemma.

\begin{lem} \label{lem.tB}
 Let $B _1, \dotsc, B _r$ be the matrices in \eqref{eq:SPj_B} and let
 $\tilde B$ be the matrix in \eqref{eq.mat-form-ps}.  Then $\mu$
 is an eigenvalue of $\tilde B$ if, and only if, $\mu ^r$ is an
 eigenvalue of $B _{r-1}B _{r} B _1\dotsm B _{r-2}$.
 In other words, the eigenvalues of $\tilde B$ are the complex roots of the eigenvalues of $B _{r-1}B _{r} B _1 \dotsm B _{r-2}$.
\end{lem}

\begin{lem}[see \S51 in \cite{Wilkinson65}] \label{lem.ABBA}
 Let $\mathcal{A}$ and $\mathcal{B}$ be square matrices. Then the spectrum of $\mathcal{AB}$ is the same as
 the spectrum of $\mathcal{BA}$.
\end{lem}
 
Combining Lemmas~\ref{lem.PPj}, \ref{lem.tB}, and \ref{lem.ABBA}, we prove Proposition~\ref{prop.BtB}. That result links the relation between the Floquet matrix with the one using the multiple shooting. To prove it, it is enough {\it i)} to observe that from Lemma~\ref{lem.PPj} the spectrum of $D _x P$ on the torus $\varphi$ has the same spectrum as $(D _x P _r)\dotsm (D _x P _1)$, respectively on $\varphi _r, \dotsc, \varphi _1$. {\it ii)} The different rotation operators $T _{-\rho /r}$ do not change the spectrum because of the Lemma~\ref{lem.ABBA}. {\it iii)} The Floquet $C _j$ do not change the spectrum of $B _j$. Therefore, Lemma~\ref{lem.tB} allows to finish the proof of the Proposition~\ref{prop.BtB}.

\begin{pro} \label{prop.BtB}
  The eigenvalues of the Floquet matrix of a multiple shooting \eqref{eq:SPj_red} with $r$ sections are the complex $r$ roots of the eigenvalues of the Floquet matrix with single shooting \eqref{eq:SP_red}.
\end{pro}

\subsection{Multiple shooting applied to invariant manifolds}

If the torus is very hyperbolic, the linear approximation to the manifold in one of the sections of the torus (single shooting) is enough to globalize
the manifold with a good level of accuracy. This is
because, as the manifold is very unstable, it is sufficient to use the unstable
direction of the torus in one of the sections, say $\varphi_1$, to grow numerically
the manifold \cite{RosalesJJC21b}. Other works, as \cite{DuartePhD}, use multiple
shooting to compute the linear approximation to the invariant manifold.
As here we are interested in a high-order approximation to these manifolds, we need to 
compute high-order derivatives of the map. Due to the strong instability of the torus,
we have to continue with the multiple shooting scheme in order to compute
the derivatives of the maps $P_j$ accurately. Therefore, we will compute
the Taylor-Fourier expansions for the torus $\varphi_j$, $j=1,\dotsc,r$.

The parametrization of the manifold, as explained in Section~\ref{sec:Param}, is done at each of the $r$ sections. 
Let $W _j$ be a formal power expansion for each $j = 1, \dotsc, r$ of the form
\begin{equation*}
 W _j (\theta , \sigma) = \sum _{k \geq 0} a _{j,k}(\theta) \sigma ^k, \qquad \theta \in \T ^d.
\end{equation*}
We denote the truncated power expansion of $W _j$ of order $m$ by $W _{j,m}$.

Let us assume, by simplicity, that $|\mu| \ne 1$ is real. Then,
applying the invariance condition of the torus to the invariant manifold leads to the equations
\begin{equation}
\label{eq.mani-inv-eq-ps}
 \begin{split}
  P _j (W _j(\theta, \sigma), \theta) &= W _{j+1} (\theta + \rho /r, \mu \sigma), \qquad j = 1, \dotsc, r-1, \\
  P _r (W _r(\theta, \sigma), \theta) &= W _1 (\theta + \rho /r, \mu \sigma). \\
 \end{split} 
\end{equation}
The zeroth order in $\sigma$ of \eqref{eq.mani-inv-eq-ps} is just the torus $\varphi _j$ in
\eqref{eq.inv-eq-ps}, that is, $a _{j,0} = \varphi _j$. The first order in $\sigma$ in \eqref{eq.mani-inv-eq-ps} has the form
\begin{equation*}
    \begin{split}
        D _x P _j (a _{j,0}(\theta), \theta) a _{j,1}(\theta) &=  a _{j+1,1} (\theta + \rho /r) \mu , \qquad j = 1, \dotsc, r-1, \\
        D _x P _r (a _{r,0}(\theta), \theta)a _{r,1}(\theta) &= a _{1,1}(\theta + \rho/r) \mu .
    \end{split}
\end{equation*}
Using the change $a _{j,1}(\theta) = C _j(\theta) v _j$ for $j = 1, \dotsc, r$ and the definition of $B _j $ in \eqref{eq:SPj_B}, we end up with
\begin{equation*}
    \begin{split}
        B _j v _{j} &=  \mu v _{j+1} , \qquad j = 1, \dotsc, r-1, \\
        B _r v _{r} &= \mu v _{1},
    \end{split}
\end{equation*}
and by the matrix-block form in Remark~\ref{rmk.mat-form-ps}, we conclude that $v = (v _r, v_1, \dotsc, v _{r-1})$ is an eigenvector of eigenvalue $\mu$ of $\tilde B$ that, by Proposition~\ref{prop.BtB}, means that $\mu$ is a $r$ root of an eigenvalue of \eqref{eq:SP_red}.

Let us now assume that we know the functions $a _{j,k}$ for $j = 1,
\dotsc, r$ and $k = 0, \dotsc, m-1$. Then, using the induction
hypothesis, for $j = 1, \dotsc, r-1$, (and similarly for $j=r$)
\begin{equation} \label{eq.mani-taylor-expa}
\begin{split}
 P _j (W _{j,{m-1}}(\theta, \sigma), \theta)
 & =P _j (W _{j,m-1}(\theta, \sigma), \theta) + 
 D _x P _j (W _{j,m-1}(\theta, \sigma), \theta) a _{j,m}(\theta) \sigma ^m + \mathcal{O}(\sigma ^{m+1}) \\ &= 
 P _j (W _{j,m-1}(\theta, \sigma), \theta)  + A _j(\theta) a _{j,m}(\theta) \sigma ^m + \mathcal{O}(\sigma ^{m+1}) \\ &= 
 W _{j+1,m-1}(\theta + \frac{\rho}{r}, \mu\sigma) + b _{j,m}(\theta ) \sigma^m+A _j(\theta) a _{j,m}(\theta)  \sigma ^m + \mathcal{O}(\sigma ^{m+1}),
\end{split}
\end{equation}
with $A _j(\theta) = D _x P _j (a _{j,0}(\theta), \theta )$.  Equating the order $\sigma ^m$ in
\eqref{eq.mani-taylor-expa} with the right hand side in \eqref{eq.mani-inv-eq-ps},
we end up with the expressions
\begin{equation*}
\begin{split}
 b _{j,m} (\theta ) + A _j (\theta) a _{j,m}(\theta) &= a _{j+1,m}(\theta + \rho /r) \mu ^m, \qquad j = 1, \dotsc, r-1, \\
 b _{r,m} (\theta ) + A _r (\theta) a _{r,m}(\theta) &=  a _{1,m}(\theta + \rho /r) \mu ^m.
\end{split} 
\end{equation*}
Introducing the Floquet change 
$a _{j, m}(\theta) = C _j(\theta) u_{j, m}(\theta)$ for all $j = 1, \dotsc, r$, 
and using
\eqref{eq:SPj_B}, we deduce
\begin{equation} \label{eq.mani-cohom3-ps}
\begin{split}
 C_{j+1}(\theta + \rho/r)^{-1} b _{j,m} (\theta ) + B _j u _{j,m}(\theta) &= \mu ^m u _{j+1,m}(\theta + \rho /r) , \qquad j = 1, \dotsc, r-1, \\
 C_{1}(\theta + \rho/r)^{-1} b _{r,m} (\theta ) + B _r u _{r,m}(\theta) &= \mu ^m u _{1,m}(\theta + \rho /r).
\end{split} 
\end{equation}
The system of equations \eqref{eq.mani-cohom3-ps} can directly be solved
as the one in \eqref{eq:tosolve}.

This scheme works for real un/stable manifolds. For the reasons discussed in Section~\ref{sec:param_stable}, to compute stable manifolds is numerically more precise to consider the inverse Poincar\'e map. We can then write similar conditions to \eqref{eq.mani-inv-eq-ps} for $P _j ^{-1}$.

\begin{rmk}
 The parametrization \eqref{eq.mani-inv-eq-ps} can also been seen as a single shooting in higher dimension. Thus we can skip the detailed
 expressions in \eqref{eq.mani-cohom3-ps} for each index $j$ and apply directly the
 Section~\ref{sec:Param} but with the matrix-block in the Remark~\ref{rmk.mat-form-ps} and with $\rho / r$ instead of $\rho$. In particular, the Algorithms~\ref{alg.mani-unstab} and \ref{alg.mani-stab} can be used with the penalty of more memory usage. Indeed, \eqref{eq.mani-cohom3-ps} can be rewritten as
 \[ \mu ^m \tilde u _m(\theta + \rho/r) = \tilde B \tilde u _m(\theta) + \tilde g _m(\theta) \]
 where $\tilde B$ is defined in Remark~\ref{rmk.mat-form-ps},
 \[
 \tilde u_m(\theta) = 
 \begin{pmatrix}
 u _{r,m}(\theta) \\ u _{1,m}(\theta) \\ \vdots \\ u _{r-1,m}(\theta)
 \end{pmatrix}, \qquad \text{and} \qquad 
 \tilde g_m(\theta) = \tilde C ^{-1}(\theta + \rho/r)
 \begin{pmatrix}
 b _{r,m}(\theta) \\ b _{1,m}(\theta) \\ \vdots \\ b _{r-1,m}(\theta)
 \end{pmatrix} .
 \]
\end{rmk}

\section{Computer implementation}\label{sec:CompImp}
This section is devoted to provide some technical details and information
concerning the computer implementation of the introduced algorithms.

First, we give some explanations about how to work with Fourier
series of several variables and the package used for this
aim. Section~\ref{subsec:JetTransport} details the idea of the jet
transport technique that we use to obtain high-order derivatives of
the Poincar\'e map through automatic differentiation. Then, since this
work is focused on the parallelism of computations,
Section~\ref{subsec:Parall} is devoted to the technical details in the
implementation of our computations and an analysis of the degree of
parallelism achieved. Finally, we include some numerical tests to
analyse the accuracy of the obtained results.

\subsection{Manipulation of Fourier series in several variables}
An effective manipulation of Fourier series has a crucial impact on the performance of Algorithms~\ref{alg.torus-floquets}, \ref{alg.mani-unstab}, and \ref{alg.mani-stab} to do the steps of shifting by $\rho$ and to solve cohomological equations. 

We do not only need to be able to express the series in its Fourier coefficients and its tabulation on a grid of angles $\theta \in \T ^d$; a process based on the Discrete Fourier Transform. We have to know how to perform operations affecting to its coefficients. In particular, we need to be able to know at each memory location which is its coefficient and vice versa. These two operations; to know the tuple $(\kappa _1, \dotsc, \kappa _d)$ from its index and to know index from its tuple, are simple and essential and they are encoded in Algorithms~\ref{alg.index2tuple-row-order} and \ref{alg.four-index2tuple}. 

The former depends on how the coefficients are packed, since it is highly dependent on the package. We are going to assume that the coefficients and its tabulation are contiguously allocated in memory but they may be interleaved and/or split in memory. Moreover, we will also assume, for simplicity, that the mesh will have odd size in each of the directions in order to avoid discussions about aliasing phenomena and Nyquist frequencies. 

\subsubsection{Truncated Fourier series representation}
Let $x$ be a real-periodic function admitting a
Fourier representation. That this,
\begin{equation}
\label{eq.fou-series}
    x (\theta) = \sum _{\kappa \in \mathbb{Z}^d} c _\kappa \exp(2 \pi \im \langle \kappa, \theta \rangle), \qquad c _\kappa \in \mathbb{C} \text{ and } \theta \in \mathbb{T}^d.
\end{equation}
where $\mathbb{T} $ is identified with $ [0, 1)$ or $[-\frac{1}{2}, \frac{1}{2})$. Note that \eqref{eq.fou-series} can easily be expressed in terms of sum of $\sin/\cos$ with real coefficients by a change of the basis.
It is standard in FFT
    (Fast Fourier Transform) algorithms to truncate the series in
    \eqref{eq.fou-series} in an even equispaced mesh $\bar N = (\bar N
    _1, \dotsc, \bar N _d) \in \mathbb{N}^d$, with $N _j = 2 \bar N _j
    + 1$ and $N = (N _ 1, \dotsc, N _d)$. Thus, the mesh points of $\theta$ are
    $\{\kappa / N \colon \kappa \in I _N\}$, being $I _N$ the set of
    indices defined by
\begin{equation}
\label{eq.set-indices}
I _N = \{ (\kappa _1, \dotsc, \kappa _d) \in \mathbb{Z}^d \colon -\bar
N _j \leq \kappa _j \leq \bar N _j \text{ for all } j = 1, \dotsc, d\}
\end{equation}
and, whose cardinal is 
\begin{equation}
\label{eq.cardIn}
    M = \prod _{j = 1} ^ d (2 \bar N _j + 1) = \prod _{j = 1} ^ d N _j.
\end{equation}
Hence, numerically \eqref{eq.fou-series} is equivalent to 
\begin{equation}
\label{eq.fou-trun-series}
    x (\theta) = \frac{1}{M} \sum _{\kappa \in I _N} \widehat x
    _\kappa \exp(2 \pi \im \langle \kappa, \theta \rangle)
\end{equation}
where $\widehat x _\kappa / M $ is close to the exact Fourier
coefficient $c _\kappa$.

Since $x(\theta)$ is always real some of the coefficients are
redundant, more explicitly $ \widehat x _{-\kappa} = \widehat x _{\kappa}
^\ast$ with ${}^\ast$ representing the conjugate as complex
numbers. This fact allows us to reduce the coefficients storage to
\begin{equation}
\label{eq.cardInast}
    \bar M = (\bar N _d + 1) \prod _{j = 1} ^ {d-1} N _j .
\end{equation}
The truncated Fourier series \eqref{eq.fou-series}, in the equispaced
mesh, is in bijection with the values of the function at those
points. That is, for $\theta \in \{\kappa / N \colon \kappa \in I _N \}$,
\begin{equation*}
    x(\theta) = \sum _{\kappa \in I _N} \widehat x _\kappa \exp(2 \pi
    \im \langle \kappa, \theta \rangle), 
    \qquad 
    \widehat x _\kappa = \frac{1}{M} \sum _{\theta \in \{\kappa / N \colon \kappa \in I _N \}} x(
    \theta) \exp(-2 \pi \im \langle \kappa, \theta \rangle).
\end{equation*}

Therefore, numerically it is equivalent to store the values $\{
x(\kappa / N) \colon \kappa \in I _N \}$ or the Fourier coefficients
$\{ \widehat x _\kappa \colon \kappa \in I _N \}$. The former contains
$M$ real numbers given in \eqref{eq.cardIn} and the latter $2 \bar M$
given in \eqref{eq.cardInast} real numbers. Note that $M$ and $2 \bar
M$ differ, and then, both representation can be stored in-place if a
padding for the table of values is considered.

\subsubsection{Major orders of Fourier coefficients}
In the literature one can mainly find two approaches to store the
Fourier coefficients depending on the memory order of them. In this work, we
have used the one followed in the FFTW3 (Fastest Fourier Transform in
the West) package \cite{FFTW05} and, in particular, in its
\texttt{dft\_c2r} and \texttt{dft\_r2c} plans.

In computer science, the ordering of data can generically be in
row-major or col-major orders. To detail the access in the Fourier
coefficients for an arbitrary dimension, and following the approach in
the FFTW3 package, we use the row-major order.

\begin{dfn}[Row-major and col-major orders]
\label{def.row-col-order}
Let $x$ be a contiguous allocated $M$-dimensional vector and let $(N
_1, \dotsc, N _d)$ be in $\mathbb{N}^d$ so that $M = \prod _{j = 1}^d
N _j$.
\begin{enumerate}
\renewcommand{\theenumi}{\roman{enumi}}
    \item $x$ is said to be in row-major order if, and only if, for
      each index $\ell = 0, \dotsc, M-1$ of $x$, there is a unique
      $(\kappa _1, \dotsc, \kappa _{d})$ in $\prod _{j = 1}^d \{0,
      \dotsc, N _j -1\}$ such that
    \begin{equation}
    \label{eq.row-order}
        \ell = (\dotsb (\kappa _1 N _2 + \kappa _2) N _3 + \dotsb) N _d +
        \kappa _d.
    \end{equation}
    \item $x$ is said to be in col-major order if, and only if, for
      each index $\ell = 0, \dotsc, M-1$ of $x$, there is a unique
      $(\kappa _1, \dotsc, \kappa _{d})$ in $\prod _{j = 1}^d \{0,
      \dotsc, N _j -1\}$ such that
    \begin{equation*}
        \ell = (\dotsb (\kappa _d N _{d-1} + \kappa _{d-1}) N _{d-2} +
        \dotsb) N _1 + \kappa _1.
    \end{equation*}
\end{enumerate}
\end{dfn}
Note that we can go from $\ell$ to its tuple $(\kappa _1, \dotsc, \kappa
_d)$ using just integer
operations. Algorithms~\ref{alg.index2tuple-row-order} details the
steps to get the unique tuple corresponding to a certain index of a
vector $x$ in row-major order (similarly for col-major order).

\begin{alg}[From index to tuple in row-major order] \
\label{alg.index2tuple-row-order}
\begin{enumerate}
\renewcommand{\theenumi}{\arabic{enumi}}
    \item [$\star$] \texttt{Input:} Integer $\ell$ and odd integers $(N
      _1, \dotsc, N _d)$.
    \item [$\star$] \texttt{Output:} Integers $(\kappa _1, \dotsc,
      \kappa _d)$ such that \eqref{eq.row-order} is verified.
    \item [$\star$] \texttt{Notation:} $/$ denotes the integer
      division and $\%$ the integer modulus.
    \item $ j \gets \ell $.
    \item For $i = d, \dotsc, 1$
    \begin{enumerate}
\renewcommand{\theenumii}{\alph{enumii}}
\renewcommand{\labelenumii}{\theenumii)}
        \item $ \kappa _i \gets j \% N _i $.
        \item $ j \gets j / N _i $.
    \end{enumerate}
\end{enumerate}
\end{alg}

Algorithms~\ref{alg.index2tuple-row-order} is suitable to recover the
tuple of mesh points from the contiguously allocated vector of them.
However, the Fourier coefficients are indexed by $I _N$ defined in
\eqref{eq.set-indices} and the
Algorithm~\ref{alg.index2tuple-row-order} does not apply directly. A
slightly modified version given in
Algorithm~\ref{alg.four-index2tuple} takes into account the different
values of the tuple $(\kappa _1, \dotsc, \kappa _d)$ in $I _N$. The
inverse process of the Algorithm~\ref{alg.four-index2tuple} is now
straightforward, in practice, this process of given an element in $I
_N$ to determine the corresponding index $\ell$ of the vector (in
row-major order) has not been needed in all the implementation of
computing the invariant torus, its Floquet, and its invariant
manifolds.

\begin{alg}[Fourier coefficients: from index to tuple] \
\label{alg.four-index2tuple}
\begin{enumerate}
\renewcommand{\theenumi}{\arabic{enumi}}
    \item [$\star$] \texttt{Input:} Integer $\ell$ and odd integers $(N
      _1, \dotsc, N _d)$.
    \item [$\star$] \texttt{Output:} Integers $(\kappa _1, \dotsc,
      \kappa _d)$ in $I _N$ defined in \eqref{eq.set-indices}.
    \item [$\star$] \texttt{Notation:} $/$ denotes the integer
      division and $\%$ the integer modulus.
    \item $ j \gets \ell $.
    \item For $i = d, \dotsc, 1$
    \begin{enumerate}
\renewcommand{\theenumii}{\alph{enumii}}
\renewcommand{\labelenumii}{\theenumii)}
        \item $ \kappa _i \gets j \% N _i $.
        \item If $\kappa _i > N _i/2$, then $\kappa _i \gets N _i - \kappa _i$.
        \item $ j \gets j / N _i $.
    \end{enumerate}
\end{enumerate}
\end{alg}
Algorithm~\ref{alg.four-index2tuple} allows us to easily apply a rigid
rotation to the quasi-periodic orbit $x \colon \mathbb{T}^d \rightarrow
\mathbb{R}$. Indeed, if $\rho$ is in $\R^d$, then
\begin{equation*}
    x(\theta + \rho) = \sum _{\kappa \in I _N} \exp (2\pi \im \langle
    \kappa, \rho \rangle) \widehat x _\kappa \exp (2\pi \im \langle
    \kappa, \theta \rangle).
\end{equation*}
and with Algorithm~\ref{alg.four-index2tuple} we can get the tuple
$\kappa$ from the index $\ell$ of the contiguously allocated vector of
size $\bar M$ encoding the Fourier coefficients. Similarly, we can
easily solve the different cohomology equations in
Algorithms~\ref{alg.torus-floquets}, \ref{alg.mani-unstab}, and \ref{alg.mani-stab}.

\subsection{Jet transport}\label{subsec:JetTransport}
Jet transport is a computational technique to obtain high-order
derivatives of the flow of an ODE with respect to initial data and/or
parameters
\cite{BerzM98,AlessiFJSV08,AlessiFJSV09,ArmellinLBB10,WilczakZ12}.  It
is based on applying automatic differentiation \cite{Griewank00} to a
numerical integrator of ODEs.  This is done by substituting the basic
arithmetic by an arithmetic of (truncated) formal power series in
several variables. A formal power series codifies the value of a
function (the constant term) and their derivatives (the coefficients
of each monomial) up to a given order (the truncation order of the
series), and the propagation of these power series through the
numerical integration produces exactly the same results as the
integration of the corresponding high-order variational equations of
the ODE (see \cite{TRAJE2019} for more details).

Here we approximate the un/stable invariant manifolds by truncated Taylor-Fourier series of $W(\theta, \sigma)$ in \eqref{eq:param} and we need to compute $P(W(\theta,\sigma), \theta)$ as a truncated Taylor-Fourier series of 
\[
P(W(\theta,\sigma), \theta)=\sum_{k\ge 0} b_k(\theta)\sigma^k,\qquad\theta\in\T^d,
\]
where $P$ is a Poincar\'e map.  The idea is to use the
equivalence between trigonometric polynomials and a suitable table of
values. As the Fourier series $a_k$ in \eqref{eq:param} are, in fact, trigonometric
polynomials, we can represent them as a suitable tabulation,
$\{a_k^{(\ell)}\}_\ell=\{a_k(\theta_\ell)\}_\ell$. Then, the Taylor-Fourier series
can be represented as a set of Taylor expansions,
\[
W^{(\ell)}(\sigma)=\sum_{k\ge 0} a_k^{(\ell)}\sigma^k,\qquad \ell=0,1,\ldots,
\]
Then, we use jet transport to compute $P(W^{(\ell)}(\sigma), \theta _\ell)$,
\[
P(W^{(\ell)}(\sigma), \theta _\ell)=\sum_{k\ge 0} b_k^{(\ell)}\sigma^k,\qquad \ell=0,1,\ldots,
\]
and, from the coefficients of these Taylor expansions, we can use
Fourier transforms to recover the corresponding Fourier series $b_k$
so that we have $P(W (\theta, \sigma),\theta)$. If, for some $b_k$, the size of the last
Fourier modes is not small, it means that more Fourier modes are
needed. Then, we restart the calculation using a finer discretization
for the corresponding $a_k$ (i.e. a larger number of points
$\theta_\ell$) so that more Fourier modes for $b_k$ are
obtained. Consequently, function $a_k$ is also approximated by a
higher number of Fourier modes.

\subsection{Parallelism}\label{subsec:Parall}
Algorithms \ref{alg.torus-floquets}, \ref{alg.mani-unstab}, and \ref{alg.mani-stab} contain steps that are highly parallelizable. Such a parallelism was already exploited in \cite{JorbaOlmedo2009} in
Algorithm~\ref{alg.torus-floquets} using the PVM library
\cite{GeistBDJMS1994} running on
a cluster of PCs connected through an Ethernet network. Here we use OpenMP~4.5 \cite{openmp45} which runs concurrently in a PC with several CPUs and it provides an easier and efficient parallelism programming.

The use of profilers for the experiments in Section~\ref{sec:Apps} shows that more than the 98\% (in both algorithms) is spent in the evaluation of the discrete map $P$ and its derivatives, which involves ODE integrations, and a lower percentage is required to solve the cohomological equations. Therefore, the parallelism strategy has consisted in running the evaluation of the
ODE integrator, Taylor \cite{JorbaZ05} and Runge-Kutta-Verner 8(9)
\cite{Verner78} in the experiments, sequentially and independently in
each of the different available CPUs of the PC. This provides an
automatic parallelism since the algorithms require to evaluate the
discrete map for each of the different angle values in the mesh in
$\mathbb{T}^d$. Note that with this approach the use of jet transport
does not provide any downside because we do not parallelize the
integrator itself.

The second level of parallelism is in the cohomological equations, the
shifting by $\rho$, and some of the matrix-solvers that are
independent to each other either in a Fourier representation or in a
table of values.

Finally, we parallelize the transformation between the table of values
and the Fourier coefficients (and vice versa). That has been done by
the feature already provided in the FFTW3 package and in combination
with the OpenMP. We did not detect a major improvement because the
package itself is already optimized enough and already the profiler
indicated that these transformations do not contribute too much in the performance when one use the FFTW3.

We took care of the potential overhead in the initialization of the
threads, that is, the different (sub)processes that are executed in
the CPUs. Thus, we initialize the threads at the beginning of the
algorithms to have ready the pool of threads and bifurcate the code
execution when we reach those parallelizable steps.

\subsection{Accuracy tests} \label{sec:Accu}

In order to ensure that the computations are correct we implement some
tests. First, regarding to the torus and the Floquet change, we
implement three tests, two of them already introduced in
\cite{OlmedoPhD}, let us call them Test~\ref{test1}, Test~\ref{test2},
and Test~\ref{test3}. Secondly, Test~\ref{test4}, is implemented to assess the 
parametrization of the hyperbolic
invariant manifolds. The Tests \ref{test2} and \ref{test3} are run after the solutions have
been obtained, which means that they can be used to check how good are
these solutions in terms of the invariance equations they must
satisfy.

In all the tests we are going to use norms and tolerances that must be chosen depending on the model, arithmetic precision, and matching with other tolerances in the algorithms, such as, the one for the Newton's process or the ODE integrations. In Section~\ref{sec:Apps}, we will made explicit all these freedoms.

\subsubsection{The invariance equation}

Algorithm~\ref{alg.torus-floquets} stops when the invariance condition \eqref{eq:InvCond} for the torus $\varphi$ and for the Floquet change \eqref{eq:B} are satisfied within a certain threshold. On the other hand, Algorithm~\ref{alg.mani-unstab} and Algorithm~\ref{alg.mani-stab} are not iterative processes and the
steps in each algorithm are deduced by imposing (by power matching in $\sigma$) the invariance equations \eqref{eq:InvEq_manitorus} and \eqref{eq:InvEq_manitorus_inverse}, respectively.

\begin{test} \label{test1}
Let $\mathcal{A}$ be a mesh in $\T ^d$. A function $z$ is said to verify the equation $\mathcal{I}(z(\theta))=0$ with tolerance
$\tau$ if, and only if,
\[
\max _{\theta\in \mathcal{A}}
  \ \|\mathcal{I}(z(\theta))\| \leq \tau.
\]
\end{test}
Note that Test~\ref{test1} can be defined in terms of the relative error instead of the absolute error.

\subsubsection{The tail of the Fourier discretization}

The test consists in checking that the truncated Fourier
representation is accurate enough with the mesh size. We use the fact that, under a smoothness assumption, the Fourier coefficients decay. In the applications in Section~\ref{sec:Apps}, these functions are analytic and then their Fourier coefficients decay exponentially.
The truncation error is approximated by the size of the last Fourier coefficients in its representation. To
prevent potential symmetries that make zero some of the entries, we
check the last two indexed coefficients.

\begin{test} \label{test2}
A truncated real Fourier representation given by
\[
x(\theta) = \frac{x^{(0)}}{2} + \sum _{|\kappa| = 1}^N x ^{(c)}_\kappa \cos \sprod{\kappa}{\theta} +  x ^{(s)}_\kappa \sin \sprod{\kappa}{\theta}, \qquad \theta \in \T ^d.
\]
is said to verify the Test \ref{test2} with
tolerance $\tau$ if, and only if, for all $\kappa \in \N^d$ such that $|\kappa | = N$ or $|\kappa | = N-1$, 
\[
 \|(x _\kappa^{(c)}, x _\kappa^{(s)})\|_2 \leq \tau.
\]
\end{test}

We apply the Test~\ref{test2} for each of the Fourier series involved in the
torus, in its Floquet change, and its parametrized manifold. Moreover,
Test~\ref{test2} can be used to keep track which of the components of the angular variables vector $\theta \in \T^d$ have the biggest
tail size in norm, and then increase the mesh size on that direction
until either we reach a maximum mesh size or we reach the desired
tolerance.

\subsubsection{The mesh}

The third test is computationally more expensive, it consists in
checking the function that we want to make zero in a different mesh 
but with the same size. A way to do this check without the need of using more computational sources is just to perform a
fixed shift by an angle, say $\gamma$, and then check if the equation
is still verified with a prescribed tolerance.
\begin{test} \label{test3}
Let $\mathcal{A}$ be a mesh in $\T ^d$. A function $z$ is said to verify the equation $\mathcal{I}(z(\theta))=0$ with tolerance
$\tau$ and shifting $\gamma  \in \T ^d$ if, and only if,
\[
\max _{\theta\in \mathcal{A}}
  \ \|\mathcal{I}(z(\theta + \gamma))\| \leq \tau.
\]
\end{test}
Note that as Test~\ref{test1}, Test~\ref{test3} can be defined in terms of the relative error instead of the absolute error.
In the case of the torus, the Test~\ref{test3} consists first in performing the
shift $\psi (\theta) = \varphi (\theta + \gamma)$ and then checking~\eqref{eq:InvCond} 
but for $\psi$ and with the same original mesh in
$\theta$, that is,
\begin{equation*}
    \max _{\theta \in \mathcal{A}} \|\psi(\theta + \rho) - P(\psi(\theta), \theta  + \gamma) \| \leq \tau.
\end{equation*}
Similarly, we can apply Test~\ref{test3} for the Floquet change $C$, and for the coefficients $a _k$ in the expansion of~\eqref{eq:param_m}.

\subsubsection{The invariant manifolds}
Once we have computed the parametrization of the manifold up to the
desired order, we check the final accuracy of the approximation of the
invariant manifold. This is done by comparing the error of the
invariance condition for the parametrization of the invariant
manifold, $W(\theta,\sigma)$, at a given angle vector $\theta$, but at
two different values of $\sigma$, say $\sigma_1$ and
$\sigma_2=\sigma_1/2$.

Notice that, when computing the parametrization of the manifold up to
order $m$, the truncation error depends on the power $m+1$ of the
parameter $\sigma$. 
\begin{test} \label{test4}
For $\sigma_i$, with $i=1,2$, we would have that
\[
\epsilon_i= | P(W_m(\theta,\sigma_i),\theta)- W_m(\theta+\rho,\lambda \sigma_i)
| \approx c\sigma_i^{m+1},
\]
where $c$ is a constant. The relation between the two errors is
\[
\frac{\epsilon_1}{\epsilon_2} \approx \frac{\sigma_1^{m+1} }{
  (\frac{\sigma_1}{2})^{m+1}} \approx 2^{m+1}.
\]
Therefore, we check that the quantity
\[
\frac{\log(\epsilon_1/\epsilon_2)}{\log(2)},
\]
has a value close to $m+1$.
\end{test}

The Test 4 will not numerically be satisfied for all values of $\sigma
_1$ and $m$ due to round-off and cancellations in $\epsilon _1$,
$\epsilon _2$, and $\epsilon _1 / \epsilon _2$. Thus, we must play
with $\sigma _1$ and $m$ in order to have enough significant digits to
avoid these digit cancellations.

\section{Applications}\label{sec:Apps}
In this section we implement two different applications whose results are provided as supplementary material. The first
example is a classical quasi-periodically forced pendulum and the
second one is an application to celestial mechanics; a model for the
Earth-Moon system subjected to five basic natural frequencies.

In order to stress the independence of the integration method, we use a Taylor integration with jet transport and tolerance $10^{-16}$ for the first method and a Runge-Kutta-Verner 8(9) with jet transport and tolerance $10^{-14}$ for the celestial mechanic one. For both of them, we use double-precision arithmetic and to verify the different tests described in Section~\ref{sec:Accu},  we consider Euclidean norms for vectors, Fr\"obenius for matrices, and a generic test tolerance of $\tau=10^{-10}$.

In all the experiments we have used the \texttt{gcc} compiler, version 8.3.0, 
on a Linux computer with two Intel(R) Xeon(R) CPU E5-2680 @2.70GHz processors, 
which give a total of 16 cores. For the sake of simplicity, in what follows we 
use the terms core and processor equally, to refer to a single 
computational unit.

\subsection{A quasi-periodically forced pendulum}\label{subsec:Pendulum}
This first application considers one of the examples included in
\cite{JorbaOlmedo2009}. The system describes the
movement of a quasi-periodically forced pendulum
\begin{equation}
\begin{aligned}
\dot{x} &=y\\
\dot{y} &=-\alpha \sin x + \varepsilon \zeta (\theta_0, \dotsc,\theta_d),\\
\dot{\theta}_i &= \omega_i, \ \  i=0,\dotsc,d
\end{aligned}
\label{eq:FP}
\end{equation}
where $x,y \in \R$, and $\alpha$ is a parameter whose value
is chosen as $0.8$. For $i=0,\dotsc,d$, $\theta_i \in \T$ and
$\varepsilon$ accounts for the weight of the forcing function $\zeta$:
\begin{equation*}
    \zeta (\theta_0,\dotsc,\theta_d)=\left[ d + 2 + \sum_{i=0}^{d}
      \cos \theta_i \right]^{-1}.
\end{equation*}
As frequencies we have chosen, with $d = 4$,
\begin{equation}
\omega_0=1, \ \ \omega_1=\sqrt{2}, \ \ \omega_2=\sqrt{3},
\ \ \omega_3=\sqrt{5}, \ \ \omega_4=\sqrt{7}.
\label{eq:FP_freq}
\end{equation}

We have applied the methodology summarized in Section~\ref{sec:RedSys}, Algorithm~\ref{alg.torus-floquets},
to obtain the torus, the Floquet transformation, and the Floquet matrix
near $x=\pi$, $y=0$ for $\varepsilon = 0.01$. 
Recalling the Section~\ref{sec:intro}, according to the dimension of the frequency vector selected, $d$, the dimension of the resulting torus of the flow~\eqref{eq:FP} near the point $(\pi,0)$ is $d+1$. By defining a returning map $P$ to the section $\theta_0=0 \mod{2\pi}$, the dimension of the torus is reduced by one. As initial seeds we used the point $(\pi, 0)$ for the torus, the identity for the Floquet transform, and the differential of $P$ at $(\pi, 0)$ for the Floquet matrix.

Each of the angles has been discretized using $N=31$ Fourier modes, that results to a total of $2N^4=1847042$ unknowns for the torus and $4 N^4=3694084$ for the Floquet change. 
Note that a direct method to compute the torus and not using the advantage of the Floquet change needs $4N^8$ memory space which is totally unfeasible.

The Algorithm~\ref{alg.torus-floquets} was run with a Newton threshold of $10^{-10}$ and the Test~\ref{test1} for the torus is satisfied with $10^{-13}$ after 3 Newton's iteration and with $10^{-12}$ for the Floquet transformation and Floquet matrix.

After the Newton convergence and the success in the Test~\ref{test1}, we apply the Test~\ref{test2} reporting the different values for each of the angular directions, that is, respectively, $10^{-10}$, $10^{-11}$, $10^{-10}$, and $10^{-11}$. The Test~\ref{test3} is also satisfied with $10^{-11}$ for the torus and $10^{-12}$ for the Floquet transformation.

The Floquet matrix, $B$, has hyperbolic real eigenvalues
$ \lambda _s = 3.625204837874207\times 10^{-3}$ and $ \lambda _u = 2.758464817115549\times
10^2$. Note that the product $\lambda _s \lambda _u$ differs from $1$ in $10^{-11}$ which matches with the requested Newton tolerance and the results of the accuracy tests as well.  
In the Table~\ref{tab:FP_d4}, we show the computational time
required for computing the torus using different number of
processors. In the same table, the speed-up factor
is included. This factor measures the relation between the time needed
for solving the system with \texttt{p} processors with respect to the time of
the linear resolution, that is, using just one processor ($\mathtt{p}=1$). Ideally, when the parallelization is performed
with \texttt{p} processors, the time should be divided by \texttt{p}. We can see in
the table that this does not happen, specially when the number of
processors increases and so the overhead in each of the processors. Some checks have been done regarding to this;
for example, disabling the Hyper-threading of the processors the computational times
remained the same. It is noteworthy that the analysis of the
profiler to our program shows that $99.76\%$ of the
computations have been parallelized. 

\begin{table}
    \centering
    {\tt
    \begin{tabular}{rrr}
p &   Total time  & speed-up\\ \hline
 1 &  8m03s & 1.00 \\
 2 &  4m10s & 1.93 \\
 4 &  2m18s & 3.50 \\
 8 &  1m14s & 6.53 \\
16 &  39s & 12.38 \\ \hline
    \end{tabular}}
    \caption{Computational time needed for computing the torus of
      system (\ref{eq:FP}) with frequency vector of dimension
      $d=4$. First column corresponds to the number of processors
      used, second one to the total time employed according to the
      number of processors, and the last one designs the
      speed-up.}
    \label{tab:FP_d4}
\end{table}

We compute the approximations to the un/stable invariant manifolds up to order 10 following the Algorithms~\ref{alg.mani-unstab} and \ref{alg.mani-stab}. Table~\ref{tab:mani_pendulum_d4} shows the required times for these computations using different number of processors and the corresponding values for the speed-up. The Test~\ref{test1} is satisfied in relative error for each of the order in $\sigma$ starting with a $10^{-14}$ at zeroth order to $10^{-11}$ at order 10. Tests~\ref{test2}, \ref{test3}, and Test~\ref{test4} have also been successful at each of the orders.

\begin{table}
    \centering
    {\tt
    \begin{tabular}{rrr||rr}
    & \multicolumn{2}{c||}{unstable} & \multicolumn{2}{c}{stable} \\ \cline{2-5}
 p & Total time & speed-up & Total time & speed-up \\ \hline
 1 &  2h33m13s & 1.00 & 2h33m31s & 1.00 \\
 2 &  1h19m08s & 1.94 & 1h28m59s & 1.73 \\
 4 &  43m12s & 3.55 & 48m23s & 3.17 \\
 8 &  22m13s & 6.90 & 25m03s & 6.13 \\
16 &  11m09s & 13.74 & 12m32s & 12.25 \\ \hline
    \end{tabular}}
    \caption{Computational time needed for the un/stable manifolds up
      to order 10 of the torus in Table~\ref{tab:FP_d4}. First column
      corresponds to the number of processors used.}
    \label{tab:mani_pendulum_d4}
\end{table}

\subsection{A quasi-periodically perturbed model for the Earth-Moon system}\label{subsec:EMSys}
G. G\'omez, J. J. Masdemont and J. M. Mondelo developed a methodology
to generate simplified Solar Systems models (SSSM) using a set of
basic frequencies, see \cite{GomezMasdemontMondelo2002,
  MondeloPhD}. The systems of equations introduced in those works
describe the motion of a massless particle subjected to
a series of time-periodic perturbations. These models are defined in such a way that if we remove all
the time-periodic dependencies present in the SSSM, the resulting models
correspond to the well-known Restricted Three-Body Problem
(RTBP), \cite{Szebehely1967}.

Among the simplified models introduced in
\cite{GomezMasdemontMondelo2002, MondeloPhD}, special attention is
paid to the Earth-Moon case, including the gravitational effect of the Sun.
For the description of this simplified model they use five
basic frequencies for the accurate characterization of the lunar
motion. The selection of these frequencies comes from
the simplified Brown theory presented in \cite{Escobal1968}.  Their
values in terms of cycles per lunar revolution (RTBP adimensional units)
are the following:

\begin{itemize}[topsep=0pt,leftmargin=12pt,itemsep=0pt]
\item mean longitude of the Moon, $\omega_1=1$,
\item mean elongation of the Moon from the Sun, $\omega_2=
  0.925195997455093$,
\item mean longitude of the lunar perigee, $\omega_3=
  8.45477852931292 \times 10^{-3}$,
\item longitude of the mean ascending node of the lunar orbit on
  the ecliptic, $\omega_4 = 4.01883841204748\times10^{-3}$,
\item Sun's mean longitude of perigee, $\omega_5 =
  3.57408131981537 \times 10^{-6}$ .
\end{itemize}

So, this model includes the perturbative effect of the solar
gravitational field, the lunar eccentricity, inclination between the
orbital plane of the Moon and the ecliptic plane, and also between the
orbital and equatorial planes.

In order to generate the model, the authors change these frequencies to a
new basis $\nu=(\nu_1,\ldots,\nu_5)$ defined as
$\nu_1 = \omega_2$,
$\nu_2 = \omega_1 - \omega_3$,
$\nu_3 = \omega_1 - \omega_2 + \omega_4$,
$\nu_4 = \omega_1 - \omega_5$, and $\nu_5 = \omega_5 - \omega_2$,
such that when the frequencies $\nu _{1}, \dotsc, \nu _i$ are
added to the unperturbed system (Earth-Moon RTBP), the simplified models SSSM$_i$ are generated for $i = 1, \dotsc, 5$, each of them subjected to $\nu _1, \dotsc, \nu _i$ perturbations.

The equations of motion for an infinitesimal particle in these models SSSM$_i$, $i = 1, \dotsc, 5$ are introduced in terms of time-dependent functions $c^i_j$, $j = 1, \dotsc, 13$,
\begin{equation}
    \left \lbrace
    \begin{aligned}
    \ddot{x} &= c^i_1 + c^i_4 \dot{x} + c^i_5 \dot{y} + c^i_7 x + c^i_8 y + c^i_9 z + c^i_{13} \frac{\partial \Omega^i}{\partial x},\\
    \ddot{y} &= c^i_2 - c^i_5 \dot{x} + c^i_4 \dot{y} + c^i_6 \dot{z} - c^i_8 x + c^i_{10} y + c^i_{11} z + c^i_{13} \frac{\partial \Omega^i}{\partial y},\\
    \ddot{z} &= c^i_3 - c^i_6 \dot{y} + c^i_4 \dot{z} + c^i_9 x - c^i_{11} y + c^i_{12} z + c^i_{13} \frac{\partial \Omega^i}{\partial z},
   \end{aligned}
    \right.
    \label{eq:MM_ecs_sun}
\end{equation}
being
\begin{equation}
\begin{split}
    \Omega^i = \frac{1-\mu}{\sqrt{(x-\mu)^2 + y^2 + z^2}} +& \frac{\mu}{\sqrt{(x-\mu+1)^2 + y^2 + z^2}}\\
    +& \frac{\mu_{S}}{\sqrt{(x-x^i_S)^2 + (y-y^i_S)^2 + (z-z^i_S)^2}},
\end{split}
    \label{eq:MM_omega_sun}
\end{equation}
where $\mu$ is the Earth-Moon mass parameter, $\mu_S$ is the mass of the Sun with respect to the sum of masses of Earth and Moon, and $x^i_S,y^i_S,z^i_S$ denote the positions of the Sun with respect to the Earth-Moon barycenter.

The quasi-periodic time-dependent functions $c^i_j$ can be computed in terms of the positions, velocities, accelerations, and over-accelerations of the two selected primaries. 
The description of these time-dependent functions as well as the positions $x^i_S$, $y^i_S$ and $z^i_S$, consists on a refined Fourier analysis, detailed in \cite{GomezMondeloFourier1, GomezMondeloFourier2}. 
The Fourier analysis results for these time-dependent functions can be found in Tables~6 and 7 of \cite{GomezMasdemontMondelo2002} and, in a more extended version, in \cite{MondeloPhD}.
Note that, regardless of the model $i$, taking $c_j=0$ except for $c_5=2$, $c_7=c_{10}=c_{13}=1$ and omitting the last term in (\ref{eq:MM_omega_sun}), the system of equations in (\ref{eq:MM_ecs_sun}) becomes that of the RTBP.

The RTBP presents five equilibrium points
(\cite{Szebehely1967}), $L _{1, \ldots, 5}$. In the Earth-Moon RTBP, $L _{1,2,3}$, are colinear of center$\times$center$\times$saddle type, and the
other two, $L _{4,5}$, form an equilateral triangle having a
dynamics of center$\times$center$\times$center. The dynamics of the
saddle parts of $L _{1,2}$ are numerically difficult to compute
because its unstable parts are of order $10^8$ and $10^6$
respectively.

The angular dimension of these points increases as the frequencies of
the SSSM are included. In the SSSM$_1$ the equilibrium points become
periodic orbits, in the SSSM$_2$ become two-dimensional quasi-periodic solutions (or 2D tori), and so
on. A way of computing these quasi-periodic solutions is to continue them from one SSSM$_i$ to SSSM$_{i+1}$ as the
number of considered frequencies increases. This continuation is
sometimes difficult due to appearance of resonances. This phenomena was studied in works like \cite{OlmedoPhD, HaroL07}.

In order to avoid the continuation problems, we add a small
dissipation parameter to the equations of the system when continuing from SSSM$_i$ to
SSSM$_{i+1}$. Thus, elliptic eigenvalues become hyperbolic and
difficulties of convergence with the algorithm coming from possible
resonances are likely removed. Once we have the invariant torus
in the system SSSM$_{i+1}$ plus the dissipation parameter, we remove
that parameter and refine the invariant object in the original SSSM$_{i+1}$.

To prevent the numerical difficulties coming from the strong instability in $L _{1,2}$, we use multiple shooting with $r$ sections, in particular, $r=4$ and $r=3$ respectively. Then we perform the computation of the torus, its
Floquet change, and Floquet matrix until we reach the SSSM$_3$
model. With this, we have obtained the invariant tori that replace $L_{1,2}$ 
in the SSSM$_3$ model, which are tori of dimensions 3 for the flow,
and their Floquet matrices. Note that, as
we are using multiple shooting, we have computed $r$ sections of the torus.
Next, we have computed the unstable manifold of each torus.
Table~\ref{tab.SSSM3} shows the computational times and corresponding
speed-up for the approximation of the unstable invariant manifolds up to order $10$.

\begin{table}[ht] 
 \centering {\tt
 \begin{tabular}{rrr||rr}
  & \multicolumn{2}{c||}{$L _1$ unstable} & \multicolumn{2}{c}{$L _2$ unstable} \\ \cline{2-5}
  $\lambda _u$ &  
  \multicolumn{2}{c||}{{\tt 1.469645480926268e+02}} & \multicolumn{2}{c}{{\tt 1.343539917760893e+02}} \\ \hline \hline
p &   Total time  &   speed-up &   Total time  &   speed-up\\ \hline
1 & 16m29s & 1.00 & 5h40m46s & 1.00 \\
2 &  8m44s & 1.89 & 2h58m09s & 1.91 \\
4 &  4m47s & 3.45 & 1h38m36s & 3.46 \\
8 &  2m30s & 6.59 & 51m39s & 6.60 \\
16 &  1m18s & 12.68 & 26m02s & 13.09 \\ \hline
  \end{tabular}}
 \caption{Computational time with \texttt{p} CPUs of the unstable manifolds
   of $L _1$ and $L _2$ of SSSM$_3$ using meshes $N=(43,43)$ and
   $N=(223,223)$, and parallel sections $4$ and $3$
   respectively.} \label{tab.SSSM3}
\end{table}

\section{Conclusions and future work}
\label{sec:conclu}
This paper has shown that the computation of high-order Taylor-Fourier expansions of un/stable invariant manifolds associated with high-dimensional tori are, nowadays, feasible. Even when the instability of the torus is very strong, we have combined the algorithms with multiple shooting methods. We have provided explicit algorithms to compute all these invariant objects.

The developed methods look efficient enough to address computation of invariant manifolds generated by several eigendirections. We plan to modify the current code for such a context as well as to manage some of the eigenvalue cases not included here. 
Similar ideas can be applied to the case when the frequency vector $\rho$ is not known or even when the internal dynamics is not a fixed rotation $\rho$, in particular, in a context when the dynamical system is autonomous. Some results in these directions have already been worked in \cite{OlmedoPhD}.

The method is highly parallelizable to compute torus, the Floquet transformation, and its invariant manifolds. In the experiments, we used OpenMP showing a really good speed-up. We are also aware of other approaches that can take advantage of the intrinsic parallelism of the algorithms such as a GPU approach. We plan  exploring in future works a GPU parallelization scheme and providing experiments showing that there is no relevant penalty in the communication between the CPU and the GPUs.

\section*{Acknowledgements}
We thank the anonymous referees for the valuable improvements in our initial manuscript.

The project leading to this application has received funding from the European Union’s Horizon 2020 research and innovation programme under the Marie Sk\l{}odowska-Curie grant agreement No 734557, the Spanish grant PGC2018-100699-B-I00 (MCIU/AEI/FEDER, UE), and the Catalan grant 2017 SGR 1374. J.G. has also been supported by the Italian grant MIUR-PRIN 20178CJA2B ``New Frontiers of Celestial Mechanics: theory and Applications'' and with funds from NextGenerationEU within the Spanish
national Recovery, Transformation and Resilience plan. B. N. is supported by the Ministry of Economy, Industry and Competitiveness of Spain through the National Plan for I+D+i (MTM2015-67724-R) and through the national scholarship
BES-2016-078722. E.O. acknowledges the 'Severo Ochoa Centre of Excellence' accreditation (CEX2019-000928-S).

\bibliographystyle{alpha}
\bibliography{refs}

\end{document}